\documentclass{article}

    \usepackage[preprint]{neurips_2024}

\usepackage{multirow}
\usepackage[cp1250]{inputenc}
\usepackage[T1]{fontenc}
\usepackage{calligra}
\usepackage{layout}

\usepackage{setspace} % Load the setspace package

\usepackage{amsmath}
\usepackage{amsthm}
\usepackage{amssymb}
\usepackage{amsfonts}
\usepackage{mathtools}

\usepackage{color}
\usepackage[dvipsnames]{xcolor}
\usepackage{graphicx}

\usepackage{verbatim}
\usepackage{caption}

\usepackage{algorithmic}
\usepackage{algorithm}

\usepackage{thmtools} 
\usepackage{thm-restate}

\usepackage{url}

\definecolor{darkscarlet}{rgb}{0.34, 0.01, 0.1}
\definecolor{yaleblue}{rgb}{0.06, 0.3, 0.57}
\definecolor{darkpowderblue}{rgb}{0.0, 0.2, 0.6}
\definecolor{midnightblue}{HTML}{0059b3}
\definecolor{noonblue}{HTML}{e5eef7}
\definecolor{chromered}{HTML}{f14233}
\definecolor{olivedrab}{HTML}{6b8e23}

\usepackage[pagebackref=true]{hyperref}
\hypersetup{
    colorlinks=true,
    citecolor=darkscarlet,
    linkcolor=darkpowderblue,
    filecolor=magenta,      
    urlcolor=yaleblue,
    menucolor=gray,
    bookmarks=true
}

\renewcommand*{\backrefalt}[4]{%
    \ifcase #1 \footnotesize{(Not cited.)}%
    \or        \footnotesize{(Cited on page~#2)}%
    \else      \footnotesize{(Cited on pages~#2)}%
    \fi}

\usepackage{dirtytalk}

\usepackage{cleveref}

\usepackage{nicefrac}
\usepackage{adjustbox}
\usepackage{threeparttable}

\usepackage{import, ifthen}

\usepackage{tcolorbox}
\usepackage{pifont}
\definecolor{mydarkred}{RGB}{192,25,25}
\definecolor{mydarkgreen}{RGB}{25,192,25}
\definecolor{mydarkblue}{RGB}{25,25,192}

\usepackage{xspace}
\newcommand{\algname}[1]{{\color{midnightblue!70!black}\small\sf#1}\xspace}

\newcommand{\norm}[1]{{\left\| #1 \right\|}}

\newcommand{\normsq}[1]{{\left\| #1 \right\|^2}} % squared norm
\newcommand{\sqn}[1]{{\left\| #1 \right\|^2}} % squared norm (alternative)
\newcommand{\inner}[2]{\left\langle #1, #2\right\rangle} % inner product

\newcommand{\sbr}[1]{\left[#1\right]} % square bracket
\newcommand{\abr}[1]{\left\langle#1\right\rangle} % angle bracket
\newcommand{\roundbr}[1]{\left(#1\right)} % round bracket
\newcommand{\curlybr}[1]{\left\{#1\right\}} % curly bracket
\newcommand{\cD}{\mathcal{D}}

\newcommand{\cF}{\mathcal{F}}
\newcommand{\cG}{\mathcal{G}}

\newcommand{\R}{\mathbb{R}}  % set of reals

\newcommand{\del}[1]{}

\newcommand{\eqdef}{:=}

\newcommand{\Exp}[1]{{\rm E} \left[ #1 \right]} 
\newcommand{\ExpCond}[2]{{\rm E}\left[\left.#1\right\vert#2\right]}
\newcommand{\ExpSub}[2]{{\rm E}_{#1}\left[#2\right]}

\newcommand{\aprox}[2]{{\rm a\text{-}prox}_{#1}\left(#2\right)}
\newcommand{\ProxSub}[2]{{\rm prox}_{#1}\left(#2\right)}

\newcommand{\pr}[1][]{
  \ifthenelse { \equal{#1}{} }
  { \ensuremath{\mathrm{P}} }
  { \ensuremath{\mathrm{P}\left(#1\right)} }
}

\theoremstyle{plain}
\newtheorem{theorem}{Theorem}\numberwithin{theorem}{section}
\numberwithin{claim}{section}
\newtheorem{remark}{Remark}\numberwithin{remark}{section}
\numberwithin{fact}{section}
 \numberwithin{exercise}{section}
\numberwithin{example}{section}
\newtheorem{lemma}[theorem]{Lemma}
\newtheorem{proposition}[theorem]{Proposition}
\newtheorem{corollary}[theorem]{Corollary}

\newtheorem{assumption}{Assumption}%\numberwithin{assumption}{section}
\theoremstyle{definition}
\newtheorem{definition}{Definition}\numberwithin{definition}{section}

\usepackage[colorinlistoftodos,bordercolor=orange,backgroundcolor=orange!20,linecolor=orange,textsize=scriptsize]{todonotes}

\usepackage{amsmath,amsfonts,bm}

\def\1{\bm{1}}

\DeclareMathAlphabet{\mathsfit}{\encodingdefault}{\sfdefault}{m}{sl}
\SetMathAlphabet{\mathsfit}{bold}{\encodingdefault}{\sfdefault}{bx}{n}

\def\sR{{\mathbb{R}}}

\DeclareMathOperator*{\argmin}{arg\,min}

\usepackage[T1]{fontenc}    % use 8-bit T1 fonts
\usepackage{hyperref}       % hyperlinks
\usepackage{url}            % simple URL typesetting
\usepackage{booktabs}       % professional-quality tables
\usepackage{amsfonts}       % blackboard math symbols
\usepackage{nicefrac}       % compact symbols for 1/2, etc.
\usepackage{microtype}      % microtypography
\usepackage{xcolor}         % colors

\usepackage{makecell}

\newcommand{\server}{{\color{black}\underline{The server:} }}
\newcommand{\client}{{\color{black}\underline{The selected client:} }}
\newcommand{\clients}{{\color{black}\underline{The selected clients:} }}

\title{SPAM: Stochastic Proximal Point Method with Momentum Variance Reduction for Non-convex Cross-Device  Federated Learning}

\newcommand{\greencheckmark}{{\color{midnightblue!70!black}\textbf{\ding{52}}}}

\newcommand{\redtimes}{{\color{red}\textbf{\ding{55}}}}

\newcommand\br[1]{\left ( #1 \right )}

\usepackage{todonotes}

\author{%
    Avetik Karagulyan%\thanks{Use footnote for providing further information about author (webpage, alternative address)---\emph{not} for acknowledging
  \qquad
  Egor Shulgin 
  \qquad
  Abdurakhmon Sadiev 
  \qquad
  Peter Richt{\'{a}}rik  \vspace{5pt} \\
  King Abdullah University of Science and Technology (KAUST) \\
  Thuwal, Saudi Arabia\\
}

\begin{document}

\maketitle

\begin{abstract}
  Cross-device training is a crucial subfield of federated learning, where the number of clients can reach into the billions. Standard approaches and local methods are prone to issues such as client drift and insensitivity to data similarities. 
  We propose a novel algorithm (\algname{SPAM}) for cross-device federated learning with non-convex losses, which solves both issues. We provide sharp analysis under second-order (Hessian) similarity, a condition satisfied by a variety of machine learning problems in practice. 
  Additionally, we extend our results to the partial participation setting, where a cohort of selected clients communicate with the server at each communication round.
  Our method is the first in its kind, that does not require the smoothness of the objective and provably benefits from clients having similar data.  
\end{abstract}

\addtocontents{toc}{\protect\setcounter{tocdepth}{1}}

\section{Introduction}

% \paragraph{Stochastic optimization and federated learning.}
Federated learning (FL) \cite{kairouz2021advances, konecny2017federated, mcmahan2017communication} is a machine learning approach where multiple entities, known as \emph{clients}, work together to solve a machine learning problem under the guidance of a \emph{central server}. 
Each client's raw data stays on their local devices and is not shared or transferred; instead, focused updates intended for immediate aggregation are used to achieve the learning goal \cite{kairouz2021advances}.
% The data of the clients is not independently or identically distributed \cite{kairouz2021advances}.

This paper focuses on \textit{cross-device} training \cite{karimireddy2021breaking}, where the clients are mobile or IoT devices. 
To model such a large number of clients, we study the following stochastic optimization problem:
\begin{equation} \label{eq:device_problem}
\min_{x \in \R^d} f(x), \quad \text{where} \quad f(x) \eqdef \ExpSub{\xi \sim \mathcal{D}}{f_\xi(x)},
\end{equation}
where $f_\xi$ may be non-convex.
Here, we do not have access to the full function $f$, nor its gradient. 
This reflects the cross-device setting, where the number of clients is extremely large (e.g., billions of mobile phones), so each client participates in the training process only a few times or maybe even once. 
Therefore, we cannot expect full participation to obtain the exact gradient. 

Instead of the exact function or gradient values, we can sample from the distribution $\cD$ and compute $f_{\xi}(x)$ and $\nabla f_{\xi}(x)$ at each point $x$. 
We assume that the gradient and the expectation are interchangeable, meaning $\ExpSub{\xi \sim \mathcal{D}}{\nabla f_{\xi}(x)} = \nabla f(x)$. 
In the context of cross-device training, $f_{\xi}$ represents the loss of client $\xi$ on its local data \cite{karimireddy2021breaking}.

The formulation \eqref{eq:device_problem} is more appropriate than the finite-sum (\textit{cross-silo}) formulation \cite{wang2021field}:
\begin{equation*}
\min_{x \in \R^d} f(x), \quad \text{where} \quad f(x) \eqdef \frac{1}{n} \sum_{i=1}^n f_i(x) ,
\end{equation*}
as the number of clients $n$ is relatively small, which is more relevant for collaborative training by organizations (e.g., medical \cite{ogier2022flamby}).

\paragraph{Communication bottleneck.}

In federated learning, broadcasting or communicating information between computing nodes, such as the current gradient vector or model state, is necessary. 
This communication often becomes the main challenge, particularly in the cross-device setting where the nodes are less powerful devices with slow network connections \cite{konecny2017federated, caldas2018expanding, kairouz2021advances}. %\cite{lin2017deep,zhang2017zipml,li2020federated}.
Two main approaches to reducing communication overhead are compression and local training. 
Communication compression uses inexact but relevant approximations of the transferred messages at each round. 
These approximations often rely on (stochastic) compression operators, which can be applied to both the gradient and the model. 
For a more detailed discussion on compression mechanisms and algorithms, see \cite{xu2020compressed, beznosikov2020biased, shulgin2022shifted}. % safaryan2022uncertainty 

\paragraph{Local training.}% and proximal point methods.}

The second technique for reducing communication overhead is to perform local training. %, which means obtaining higher-quality gradient information on the clients. 
Local \algname{SGD} steps have been a crucial component of practical federated training algorithms since the inception of the field, demonstrating strong empirical performance by improving communication efficiency \cite{mangasarian1993backpropagation, mcdonald2010distributed, mcmahan2017communication}. 
However, rigorous theoretical explanations for this phenomenon were lacking until the recent introduction of the \algname{ProxSkip} method \cite{mishchenko2022proxskip}. 
\algname{ScaffNew} (\algname{ProxSkip} specialized for the distributed setting) has been shown to provide accelerated communication complexity in the convex setting. 
While \algname{ScaffNew} works for any level of heterogeneity, it does not benefit from the similarity between clients. 
In addition, methods like \algname{ScaffNew}, designed to fix the client drift issue \cite{acar2020federated, karimireddy2020scaffold}, require each client to maintain state (control variate), which is incompatible with cross-device FL \cite{reddi2020adaptive}.

\paragraph{Partial participation.}

In generic (cross-silo) federated learning, periodically, all clients may be active in a single communication round.
However, an important property of cross-device learning is the impracticality of accessing all clients simultaneously. Most clients might be available only once during the entire training process. Therefore, it is crucial to design federated learning methods where only a small cohort of devices participates in each round. Modeling the problem according to \eqref{eq:device_problem} naturally avoids the possibility of engaging all clients at once. We refer the reader to \cite{reddi2020adaptive, karimireddy2021breaking, khaled2022faster} %condat2022murana, khaled2022faster,mishchenko2022proxskip, patel2022towards} 
for more details on partial participation.

\paragraph{Data heterogeneity.}  

Despite recent progress in federated learning, handling data heterogeneity across clients remains a significant challenge \cite{kairouz2021advances}. 
Empirical observations show that clients' labels for similar inputs can vary significantly \cite{arivazhagan2019federated, silva2022fedembed}. 
This variation arises from clients having different preferences. 
When local steps are used in this context, clients tend to overfit their own data, a phenomenon known as client drift.

An alternative to local gradient steps is a local proximal point operator oracle, which involves solving a regularized local optimization problem on the selected client(s). 
This approach underlies \algname{FedProx} \cite{li2020federated}, which relies on a restrictive heterogeneity assumption. 
The algorithm was analyzed from the perspective of the Stochastic Proximal Point Method (\algname{SPPM}) in \cite{yuan2022convergence}. 
Independently, the theory of \algname{SPPM} has been shown to be compatible with the second-order similarity condition (\Cref{as:similarity}) from an analytical perspective. 
Based on these connections, various studies have explored \algname{SPPM}-based federated learning algorithms, and we refer interested readers to \cite{khaled2022faster,lin2024stochastic} for more details.

\subsection{Prior work}

\paragraph{Momentum.}

Momentum Variance Reduction (\algname{MVR}) was introduced in the context of server-only stochastic non-convex optimization \cite{cutkosky2019momentum}. The primary motivation behind this method, also known as \algname{STORM}, was to avoid computing full gradients (which is impractical in the stochastic setting) or requiring "giant batch sizes" of order $\mathcal{O}(1/\varepsilon^2)$. Such large batch sizes are necessary for other methods like \algname{PAGE} \cite{li2021page} to find an $\varepsilon$-stationary point.

The authors assume bounded variance for stochastic gradients $\nabla f_{\xi}$ and analyze the method under additional restrictive conditions. 
However, these conditions can be replaced with the second-order heterogeneity (\Cref{as:similarity}).
The convergence result of \algname{MVR} for non-convex objectives includes the stochastic gradient noise term $\sigma^2$ in the upper bound. 
To eliminate the dependence on this parameter, they propose an adaptive stepsize schedule under the additional assumption that $f_{\xi}$ is Lipschitz continuous.

\paragraph{MIME.} 

\,\algname{MIME} is a flexible framework that makes existing optimization algorithms applicable in the distributed setting \cite{karimireddy2021breaking}. They describe a general scheme that combines local \algname{SGD} updates with a generic server optimization algorithm. The authors then study particular instances of the framework, such as \algname{MIME + ADAM} \cite{kingma2014adam} and \algname{MIME + MVR} \cite{cutkosky2019momentum}.

However, their analysis with local steps is limited from the non-convex cross-device learning perspective. 
First, they assume smoothness also in the case of one sampled client. 
Moreover, \algname{MIME} suffers from a common issue of local methods. 
In Theorem 4 of \cite{karimireddy2021breaking},  the stepsize is taken to be of order $O(\nicefrac{1}{Lm})$, where $L$ is the smoothness parameter of the client loss and $m$ is the number of local steps. 
This means that the stepsize is so small that multiple steps become equivalent to a single, smoother stochastic gradient descent step, negating the potential benefits of local \algname{SGD}. 
Finally, their analysis requires an additional weak convexity assumption for the objective in the partial participation setting.

\paragraph{CE-LSGD.} The Communication Efficient Local Stochastic Gradient Descent (\algname{CE-LSGD}) was introduced by \cite{patel2022towards}. 
They propose and analyze two algorithms, with the second one tailored for the cross-device setting \eqref{eq:device_problem}. 
This algorithm comprises two components: the \algname{MVR} update on the server and \algname{SARAH} local steps on the selected client. 
The latter, known as the Stochastic Recursive Gradient Algorithm, is a variance-reduced version of \algname{SGD} that periodically requires the gradient of the objective function \cite{nguyen2017sarah}. 
% Subsequently, \cite{beznosikov2024random} reduced the necessity for full gradient computation by employing a randomized reshuffling strategy and aggregating stochastic gradients obtained in each epoch.

The analysis of \cite{patel2022towards} explicitly describes how to choose the number of local updates and the local stepsize.  
They also provide lower bounds for two-point first-order oracle-based federated learning algorithms.
The drawback of their setting is that to have meaningful local updates, they need smoothness of each client function 
$f_{\xi}$. 
In addition, similar to \algname{MIME}, it has a dependence of the stepsize on the number of local steps, which undermines the benefits of doing many steps.

\paragraph{SABER.}
The \algname{SABER} algorithm combines \algname{SPPM} updates on the clients with \algname{PAGE} updates on the server \cite{mishchenko2023federated}. 
Their paper utilizes Hessian similarity (\Cref{as:similarity}) and leverages it for the finite-sum optimization objective.  
However, their analysis for the partial participation setting relies on an assumption difficult to verify in the general non-convex regime. 
In fact, if the function is not weakly convex, as in the case of \algname{MIME}, this assumption may not hold. 
Specifically, it requires that $ f\left(\frac{1}{B} \sum_{i=1}^{B} w_i\right) \leq \frac{1}{B} \sum_{i=1}^{B} f(w_i) $, where $ w_i $ are arbitrary vectors in $\mathbb{R}^d$ obtained using proximal point operators.

\subsection{Contributions}

This paper introduces a novel method called Stochastic Proximal point And Momentum (\algname{SPAM}). 
Our method combines Momentum Variance Reduction (\algname{MVR}) on the server side to leverage its efficiency in stochastic optimization while employing Stochastic Proximal Point Method (\algname{SPPM}) updates on the clients' side.
We analyze four versions of the proposed algorithm:
\begin{itemize}
    \item \algname{SPAM}  - using exact \algname{PPM} with constant parameters,
    \item \algname{SPAM} - employing exact \algname{PPM} with varying parameters,
    \item \algname{SPAM-inexact} - employing inexact \algname{PPM} with varying parameters,
    \item \algname{SPAM-PP} - using inexact \algname{PPM} with varying parameters and partial participation.
\end{itemize}

We then carry out an in-depth theoretical analysis of the proposed methods, showcasing their advantages compared to relevant competitors and addressing the limitations present in those works. 
Specifically, we demonstrate convergence upper bounds on the average expected gradient norm for all variants of \algname{SPAM}.
% The analysis includes the stationarity guarantees for all variants of \algname{SPAM}. 
% Specifically, we demonstrate the convergence of the average expected gradient norm to a neighborhood of $0$ for all variants. 
% This is the standard convergence type in non-convex optimization \cite{cutkosky2019momentum}.

We also conduct a communication complexity analysis based on our convergence results. 
Namely, we show that \algname{SPAM} can provably benefit from similarity and significantly improve upon the lower bound for centralized (server-only) methods.
% Specifically, we match the lower bounds, established in \cite{patel2022towards}, for the number of iterations required to reach precision error $\varepsilon$ for \algname{SPAM-PP}.
In addition, we design a varying stepsize schedule that removes the neighborhood from the stationarity bounds. 
% In addition, following the varying stepsize scheme introduced in the original \algname{MVR} paper, we design a stepsize schedule that removes the neighborhood from the stationarity bounds. 
Leveraging this scheme, our algorithm achieves the optimal convergence rate of $O(1/K^{1/3})$, where $K$ denotes the number of iterations.

Our algorithms, in particular \algname{SPAM-PP}, shine in the cross-device setting when compared to the competitors.
First, in contrast to non-\algname{SPPM}-based algorithms, such as \algname{MIME} and \algname{CE-LSGD}, we allow greater {\it flexibility for the local solvers}. 
Thus, unlike  \algname{MIME} and \algname{CE-LSGD}, we do not require either convexity or smoothness of the local objectives.
Our algorithm is compatible with any local solver as soon as it satisfies certain conditions outlined in \Cref{def:inexact}. 
Furthermore, compared to \algname{SABER}, our partial participation setting does not require (weak) convexity of the objective.
Moreover, we offer analysis that is substantially simpler than prior works and can be of independent interest outside of the FL context.
We present a visual comparison of the relevant methods in \Cref{table:overview}.

Another important aspect of our algorithms is that they do not need local states/control variates to be stored on each client, as opposed to many standard federated learning techniques \cite{karimireddy2020scaffold, acar2020federated, mishchenko2022proxskip}.
This is crucial for cross-device learning as each client may participate in training a single time. 
% The absence of these states allows to significantly improve the memory requirements for the algorithm.

Finally, we validate our theoretical findings through meticulously designed experiments. Specifically, we tackle a federated ridge regression problem, where we can precisely control the second-order heterogeneity parameter $\delta$, as well as the computation of the local proximal operator.

\paragraph{Paper Organization.}
The rest of the paper is organized as follows. 
\Cref{sec:assumptions} presents the mathematical notation and the theoretical assumptions we use in the analysis.
The main algorithm with the exact proximal point oracle and its theoretical analysis is presented in \Cref{sec:spam}.
The algorithm with an inexact proximal operator follows in \Cref{sec:inexact}. 
In \Cref{sec:spam-pp-new}, we show our most general algorithm, \algname{SPAM-PP}, which uses random cohorts of clients. 
We present experiments in \Cref{sec:experiments} and conclude the paper in \Cref{sec:conclusion}.

\begin{table}[t]
\addtolength{\tabcolsep}{-0.2em}
\centering
\caption{Comparison of the proposed algorithm with other relevant methods. 
%The columns are: {\bf 2ODH} - 2nd-order Data Heterogeneity, {\bf NCT} - Non-Convex Theory, {\bf PP} - Partial Participation, {\bf NSA} - No Smoothness Assumption, {\bf CD} - Cross Device, {\bf SU} - Server Update, {\bf CO} - Client Oracle.
} 
\renewcommand{\arraystretch}{1.4} 
\begin{tabular}{c | c c c c c c c}
\hline\hline
% \textbf{Algorithm} & \textbf{2ODH}  & \textbf{NCT} & \textbf{PP}   & {\bf NSA}  & {\bf CD} & {\bf SU} & {\bf CO} &  \textbf{Paper} \\\hline
\textbf{Algorithm} & 
\makecell{\textbf{Hessian} \\ \textbf{similarity}} & 
\makecell{\textbf{Partial}\\ \textbf{Participation}} & 
\makecell{\textbf{No Smoothness}\\ \textbf{assumption}} & 
\makecell{\textbf{Cross}\\ \textbf{Device}} & 
\makecell{{\bf Server} \\ \textbf{update}} & 
\makecell{{\bf Client} \\ \textbf{oracle}}\\\hline
\algname{FedProx} \cite{yuan2022convergence} & \redtimes & \greencheckmark   & \greencheckmark   & \greencheckmark & {\bf --} & \algname{PPM}  \\ 
\algname{SABER} \cite{mishchenko2023federated} & \greencheckmark & \redtimes  & \greencheckmark & \redtimes & \algname{PAGE} & \algname{PPM} \\ 
\algname{MIME} \cite{karimireddy2021breaking} & \greencheckmark & \redtimes  & \redtimes & \greencheckmark & \algname{MVR} & \algname{SGD} \\ 
\algname{CE-LSGD} \cite{patel2022towards} & \greencheckmark & \greencheckmark   & \redtimes & \greencheckmark & \algname{MVR} & \algname{SARAH} \\ 
\hline
\algname{SPAM} & \greencheckmark & \greencheckmark    & \greencheckmark & \greencheckmark & \algname{MVR} & \algname{PPM} \\ 
\hline\hline
\end{tabular}
\label{table:overview}
\end{table}

\section{Notation and assumptions}\label{sec:assumptions}

We use $\nabla f$ for the gradient, $\norm{\cdot}$ for the Euclidean norm, and $\Exp{\cdot}$ for the expectation. ${\sf Unif}(S)$ denotes uniform distribution over the discrete set $S$.
The proximal point operator of a real-valued function $g:\R^d\rightarrow \R$ is defined as the solution of the following optimization 
\begin{equation}
   \ProxSub{g}{x} := \arg\min_{y \in \R^d}\curlybr{g(y) + \frac12\normsq{x-y}}. 
\end{equation}
% $\ProxSub{g}{x} := \arg\min_{y}\curlybr{g(y) + \frac12\normsq{x-y}}$. 
We refer the reader to \cite{beck2017first} for the properties of the proximal point operator.
There exists a lower bound for function $f$, and it is denoted as $f_{\inf} > -\infty$.

We use index $i$ for a non-random client, while $\xi$ is used for a randomly selected client.
One of the main assumptions of our analysis is that we have access to stochastic samples $\xi \sim \mathcal{D}$ and in particular, we can evaluate the gradient $\nabla f_{\xi}$ at any point $x \in \R^d$. 
\begin{assumption}[\textbf{Bounded variance}]
\label{as:sigma}
    We assume there exists $\sigma \geq 0$ such that for any $x \in \R^d$
    \begin{equation}
        \label{eq:sigma}
        \Exp{\sqn{\nabla f_{\xi}(x) - \nabla f(x)}} \leq \sigma^2
    \end{equation}
\end{assumption}

We say that the function $f$ is $L$-smooth, if its gradient is Lipschitz continuous $\forall x,y\in \R^d$:
\begin{equation}\label{eq:smoothness}
    \norm{ \nabla f(y) - \nabla f(x) } \leq {L}\norm{x-y}.%, \quad \text{for} \quad \forall x,y\in \R^d.
\end{equation}
In many machine learning scenarios, the non-convex objective functions do not satisfy \eqref{eq:smoothness}. 
Moreover, several prior works \cite{zhang2019gradient, crawshaw2022robustness} showed that such smoothness condition does not capture the properties of popular models like LSTM, Recurrent Neural Networks, and Transformers.
% There is a vast body of work that studies alternative assumptions and techniques for minimization which may be more suitable for modern models \cite{zhang2019gradient, li2024convex}. 

Our second assumption is the second-order heterogeneity. 
Further in the analysis, this assumption will take the role of smoothness.
\begin{assumption}[\textbf{Hessian similarity}]
\label{as:similarity}
    Assume 
    % for each $x, y \in \R^d$ 
    there exists $\delta \geq 0$ such that for any $\xi$ and $x, y \in \R^d$
    \begin{equation}
        \label{eq:similarity_almost_sure}
        \norm{\nabla f_{\xi}(x) - \nabla f(x)  - \nabla f_{\xi}(y) + \nabla f(y)} \leq \delta \norm{x - y}.
    \end{equation}
\end{assumption}

When all functions $f_\xi$ are twice-differentiable condition \eqref{eq:similarity_almost_sure} can also be formulated as
\begin{equation} \label{eq:Hessian_similarity}
    \norm{\nabla^2 f_{\xi}(x) - \nabla^2 f(x)} \leq \delta,
\end{equation}
motivating the name \textit{second-order heterogeneity} used interchangeably with \textit{Hessian similarity} \cite{khaled2022faster}.

% \esnote{ mention first-order heterogeneity}
% \esnote{We have to explain why it is better than standard smoothness as it is more general and may better capture real-world practical scenarios as clients tend to have some level of similarity. 
% In addition, we can criticize the prior (using only $L$-smoothness) works due to their limited applicability for neural networks (transformers) as there is a better $(L_0, L_1)$ assumption \cite{zhang2019gradient}.}

This assumption \cite{mairal2013optimization, shamir2014communication, mairal2015incremental} holds for a large class of machine learning problems. % where the input data are similar, but the labels vary. 
Typical examples include least squares regression,
classification with logistic loss \cite{woodworth2023two},  statistical learning for quadratics \cite{shamir2014communication}, generalized linear models \cite{hendrikx2020statistically}, and semi-supervised learning \cite{chayti2022optimization}.
Furthermore, a similar assumption was used to improve convergence results in centralized \cite{tyurin2023sharper} and communication-constrained distributed settings \cite{szlendak2022permutation}.
In the distributed setting, \eqref{eq:similarity_almost_sure} is especially relevant as the parameter $\delta$ remains small, even if different clients have similar input distributions but widely varying outputs for the same input.
% it allows to leverage similarity 
% Specifically, the parameter $\delta$ remains small, even if different clients have similar input distributions but widely varying outputs for the same input.  
See more details on the assumption in \cite[Section 9]{khaled2022faster} and \cite[Section 3]{woodworth2023two} for discussion on synthetic data, private learning, etc.% beznosikov2024similarity
% is pertinent to the field of stochastic optimization, where it is required for detailed theoretical analysis \cite{khaled2022faster,tyurin2023sharper}.
% For deterministic convex setting \cite{kovalev2022optimal} referred to \eqref{eq:Hessian_similarity} as \textit{function similarity} and obtained an optimal method without client sampling.

In the following sections, we present our main algorithms as well as the corresponding convergence theorems. We focus on the non-convex optimization problem \eqref{eq:device_problem}, where the goal is to find an $\varepsilon$-approximate stationary point $x \in \mathbb{R}^d$ such that $\Exp{\sqn{\nabla f(x)}} \leq \varepsilon$.
% We are interested in finding an approximately stationary point of the nonconvex problem

\section{SPAM}\label{sec:spam}

In this section, we describe our main algorithm in its simpler form, that is, \algname{SPAM} with one sampled client and exact proximal point computations. 
We then provide theoretical convergence guarantees and a complexity analysis of the proposed methods.

The algorithm proceeds as follows. 
We first choose a stepsize sequence $\gamma_k$ and a momentum sequence $p_k$. The server samples a client. 
The selected client then computes the new gradient estimator $g_k$ and assigns the new iterate as the proximal point operator with a shifted gradient term:
\begin{align*}
    x_{k+1} 
    &= \ProxSub{\gamma_k f_{\xi_k}}{x_k+\gamma_k(\nabla f_{\xi_k}(x_{k})-g_k)} 
    = \arg\min_{y} \phi_{k} (y),
\end{align*}
where $\phi_k$ is defined as 
\begin{equation}\label{eq:inexact}
    \phi_k (y) := f_{k}(y)+\inner{g_k-\nabla f_{k} (x_k)}{y - x_k} + \frac{1}{2\gamma_k}\|y-x_k\|^2.
\end{equation} 
 The new iterate is then sent to the server, and the process repeats itself. For the algorithm's pseudocode, please refer to \Cref{alg:spam-exact}.
\begin{algorithm}
    \begin{algorithmic}[1]
        \STATE  \textbf{Input:}  Starting point $x_0 = x_{-1} \in\R^d$,  initialize $g_0 = g_{-1}$,\\
        choose $\gamma_k > 0$ and $p_k > 0$;
        \FOR {$k=0,1,2, \ldots$}
        \STATE \server samples $\xi_k\sim \cD$; 
        \STATE \client sets $g_k =  \nabla f_{\xi_k} (x_k) +(1-p_k)\left(g_{k-1} -\nabla f_{\xi_k} (x_{k-1}) \right)$;
        \STATE \client sets $x_{k+1} \in
        \begin{cases}
         \ProxSub{\gamma_k f_{\xi_k}}{x_k+\gamma_k(\nabla f_{\xi_k}(x_{k})-g_k)}; 
            \hspace{6pt}\rhd \ \text{\hfill{\algname{SPAM}}} \\
            \aprox{\epsilon}{x_k,g_k,\gamma_k,\xi_k}; \hspace{65pt}\rhd \ \text{\algname{SPAM-inexact}}
        \end{cases}
        $
        \STATE \client sends $x_{k+1}$ to the server.
        \ENDFOR
    \end{algorithmic}
    \caption{\algname{SPAM}, \algname{SPAM-inexact}}
    \label{alg:spam-exact}    
\end{algorithm}

The following proposition is the cornerstone of our analysis. 
It provides a recurrent bound for a certain sequence $V_k$, which serves as a Lyapunov function:
\begin{equation}\label{eq:lyapunov}
    V_{k} = f(x_k) - f_{\inf} + \frac{15\gamma_k}{16(2p_k - p_k^2)}\sqn{g_k - \nabla f(x_k)}.
\end{equation}

\begin{proposition}\label{prop:lyapunov}
    Let $x_k$ be the iterates of \algname{SPAM} for an objective function $f$, which satisfies Assumptions  \ref{as:sigma} and \ref{as:similarity}. 
    If $\gamma_k^2 \leq \min\left\{\frac{1}{16\delta^2},\frac{p_k}{96\delta^2(1-p_k)}\right\}$, then for every $k \geq 1$ %\in \mathbb{N}$
    \begin{eqnarray*}
     \Exp{V_{k+1}}    
     &\leq &  \Exp{V_k} - \frac{\gamma_k}{32}\Exp{\sqn{\nabla f(x_{k+1})}} + 2\gamma_k p_k\sigma^2,
 \end{eqnarray*}
 where $V_k$ is defined in \eqref{eq:lyapunov}. 
\end{proposition}
The proof can be found in \Cref{app:proof-lyapunov}. 
This proposition leads to a convergence result for \algname{SPAM} with fixed parameters.
\begin{theorem}[\algname{SPAM} with constant parameters]\label{thm:sppm-mvr}
 Suppose Assumptions  \ref{as:sigma}, \ref{as:similarity} are satisfied. Then,
\begin{eqnarray}\label{eq:spam-thm}
    \frac{1}{K}\sum_{k=1}^{K}\Exp{\sqn{\nabla f(x_{k})}} 
    &\leq& \frac{32(f(x_0) - f_{\inf})}{\gamma K} + \frac{32\sqn{g_0 - \nabla f(x_0)}}{ (2p - p^2)K} + 64 p\sigma^2,
\end{eqnarray}
where $\gamma^2 \leq \min\left\{\frac{1}{16\delta^2},\frac{p}{96\delta^2(1-p)}\right\}$.
\end{theorem}
The proof of the theorem can be found in \Cref{app:proof-thm:sppm-mvr}.
\begin{corollary}
The result can also be written as 
    \begin{eqnarray*}
    \Exp{\sqn{\nabla f(\tilde{x}_{K+1})}} 
    &\leq& \frac{32(f(x_0) - f_{\inf})}{\gamma K} + \frac{32\sqn{g_0 - \nabla f(x_0)}}{ (2p - p^2)K} + 64 p\sigma^2,
\end{eqnarray*}
where $\tilde{x}_{K+1}$ is taken uniformly randomly from the iterates of the algorithm $\{x_1,x_2,\ldots,x_{K+1} \}$.
\end{corollary}
Our primary focus is communication complexity, which is typically the main bottleneck in cross-device federated settings \cite{kairouz2021advances}. Below, we present the communication complexity of \algname{SPAM} with fixed parameters.
\begin{corollary}\label{cor:spam-complexity}
    Choose constant stepsize $\gamma_k = \gamma$ and momentum parameter $p_k = p$ as
    \begin{equation*}
        \gamma = \min \br{\frac{1}{\delta}, \br{\frac{F}{2 \delta^ 2\sigma^2 K}}^{1/3}},
        \qquad p = \max(\gamma^2\delta^2, 1/K).
    \end{equation*}
    % $\gamma_k = \gamma$, a momentum parameter $p_k = p = \max(\gamma^2\delta^2, 1/K)$.
    % Let us fix the parameters $\gamma_k = \gamma$ and $p_k = p$, such that $\gamma$ and $p$ satisfy the conditions of \Cref{thm:sppm-mvr}. 
    Then, the communication complexity of \algname{SPAM}, to obtain $\varepsilon$ error is of order 
    $\mathcal{O} \br{\frac{\delta F + \sigma^2}{\varepsilon} + \frac{\delta \sigma F}{\varepsilon^{3/2}}}$,
    where $F \eqdef f(x_0) - f_{\inf}$.
\end{corollary}

The proof is deferred to \Cref{app:cor:spam-complexity}. Our result indicates that higher similarity (smaller $\delta$) leads to fewer communication rounds to solve the problem.
Obtained complexity remarkably improves upon the lower bound %(ignoring $F = f(x_0) - f_{\inf}$)
$\mathcal{O} \br{\frac{L F + \sigma^2}{\varepsilon} + \frac{L \sigma F}{\varepsilon^{3/2}}}$
for $\delta \ll L$ \cite{arjevani2023lower}. 

Suppose now that we can initialize $g_0 = \nabla f(x_0)$. 
% have access to the full gradient at the starting point \(x_0\). 
Then, the second term in the convergence upper bound \eqref{eq:spam-thm} vanishes.
Repeating the exact steps as in the proof of \Cref{cor:spam-complexity}, we obtain the convergence rate:
$\mathcal{O} \left(\frac{\delta F}{K} + \left(\frac{\delta \sigma F}{K}\right)^{2/3}\right)$,
which leads to a communication complexity of
$\mathcal{O} \left(\frac{\delta F}{\varepsilon} + \frac{\delta \sigma F}{\varepsilon^{3/2}}\right)$.
Thus, our result shows that in the homogeneous case (i.e., $\delta = 0$), communication is not needed at all, as each client can solve the problem locally.

\begin{remark}
    Lower bounds for two-point first-order oracle federated learning algorithms with local steps were investigated in \cite{patel2022towards}. 
    However, these bounds are specifically designed for local \algname{SGD}-type methods, such as \algname{MIME}. In addition, results by \cite{patel2022towards} require smoothness.
    As our methods are agnostic to the choice of local solvers, the applicability of these bounds to our setting remains limited.
\end{remark}

It is important to highlight that the stepsize $\gamma$ in \algname{SPAM} differs from the stepsize used in local methods such as \algname{MIME} and \algname{CE-LSGD}. In these methods, the stepsize is intended to run the algorithms locally on a selected client. However, \algname{SPAM} only requires an oracle for proximal points, allowing the oracle to use any optimization method suitable for the problem at hand. Additionally, the stepsize for local \algname{SGD}-based methods depends on the smoothness parameter, which is not required in our theorem. Thus, our approach allows much more flexibility for choosing local solvers that are adaptive to the curvature of the loss \cite{malitsky2020adaptive, mishkin2024directional}.
See \Cref{table:overview} for a detailed comparison of the methods.

 % The resulting iteration complexity of the algorithm is 
 %    \begin{equation}\label{eq:spam-complexity}
 %        K = O\roundbr{\max\curlybr{{\delta},
 %        \delta\sqrt{\frac{\sigma^2- {\varepsilon} }{\varepsilon}}, \frac{{\sigma^2}\sqn{g_0 - \nabla f(x_0)}}{{\varepsilon}}}\frac{1}{\varepsilon} }. 
 %    \end{equation}\avetik{Here we do not "really" benefit from small $\delta$, as the third element in \eqref{eq:spam-complexity} does not depend on it. }.
 %    In the case that we have access to the full gradient at the beginning of the training, then we can pick $g_0 = g_{-1} = \nabla f(x_0)$. This will eliminate the second term of \eqref{eq:spam-thm}.
 %     If $\delta$ is equal to zero, then $f_{\xi} = f$ for any $\xi$. Thus, we need only one local proximal point computation to reach the minimum.
 %    \avetik{might be wrong}
 %    .

In \eqref{eq:spam-thm}, we notice that the last term, which is due to the stochastic nature of our problem, does not vanish when $K$ is large. 
To remove the stationarity neighborhood, let us now consider varying stepsizes for \algname{SPAM}, with decaying momentum parameters $p_k$. 
\begin{theorem}[\algname{SPAM}]\label{thm:spam-decaying}
    Consider \algname{SPAM} for an objective function $f$ that satisfies Assumptions \ref{as:sigma} and \ref{as:similarity}. 
    Let $\gamma_k$ be a sequence of varying stepsizes satisfying $\gamma_k^2 \leq \frac{1}{16\delta^2}$ and choose $p_k = \frac{96\delta^2\gamma_k^2}{96\delta^2\gamma_k^2 + 1}$. 
    Then, 
    \begin{eqnarray}\label{eq:spam-decreasing-thm}
    \frac{1}{\Gamma_K}\sum_{k=1}^{K} \gamma_k\Exp{\sqn{\nabla f(x_{k})}} 
    &\leq& \frac{32V_0}{\Gamma_K} 
    + \frac{2}{\Gamma_K}\sum_{k=1}^{K} \frac{96\delta^2\gamma_k^3}{96\delta^2\gamma_k^2 + 1} \sigma^2,
\end{eqnarray}
where $\Gamma_K = \sum_{k = 1}^{K} \gamma_k$. 
\end{theorem}
The proof of \Cref{thm:spam-decaying} can be found in \Cref{app:proof-decaying}. 
\begin{remark}
Similar to \Cref{thm:sppm-mvr}, we can represent the left-hand side of \eqref{eq:spam-decreasing-thm} with a single expectation: 
$\Exp{\sqn{\nabla f(\tilde{x}_{K})}}$, where $\tilde{x}_{K} = {x}_{i}$, for $i = 1,\ldots,K$ with probability $\gamma_i/\Gamma_K$.
\end{remark}
To ensure that the right-hand side converges to zero as $ K \rightarrow \infty $, we need $\gamma_K \rightarrow 0$ and $\Gamma_K \rightarrow +\infty$. This suggests using a stepsize schedule of order $\gamma_k = O(k^{\beta-1})$, implying $\Gamma_K = O(K^{\beta})$ for some $\beta \in (0,1)$. Consequently, the right-hand side of \eqref{eq:spam-decreasing-thm} is of order $O(K^{-\beta} + K^{2\beta-2})$. By optimizing over $\beta$, we deduce that $\gamma_k = O(k^{-\nicefrac{1}{3}})$ results in a stationarity bound of order $O(K^{-\nicefrac{2}{3}})$.

\begin{corollary}[Optimal stepsize schedule]\label{cor:varying-stepsize}
    If  $\gamma_k = \frac{1}{4\delta k^{\nicefrac{1}{3}}}$ and $p_k = \frac{96\delta^2\gamma_k^2}{96\delta^2\gamma_k^2 + 1}$, then to obtain $\varepsilon$-stationarity for \algname{SPAM} 
    we need $K = O({\varepsilon}^{3/2})$ iterations under assumptions \ref{as:sigma} and \ref{as:similarity}. 
    %Here, $\delta$ and $\sigma^2$ are the parameters of assumptions \ref{as:sigma} and \ref{as:similarity}.
\end{corollary}

This coincides with the existing lower bounds by \cite{arjevani2023lower}, meaning that our result is tight up to constants.

% \subsection{Complexity analysis} \label{sec:SPAM_complexity}

% \esnote{Highlight the benefit from greater similarity (smaller $\delta$).}

% \esnote{Contrast to PAGE \cite{li2021page, tyurin2023sharper}, SABER \cite{mishchenko2023federated}, and SVRS \cite{lin2024stochastic} in details.} 

\section{Inexact proximal operator}\label{sec:inexact}

In the previous theorems, we assume that each sampled client $\xi_k$ can exactly compute the proximal operator to obtain the new iterate $x_{k+1}$.
The latter means that this client can exactly solve a (potentially) non-convex minimization problem, which might be problematic in practice.
However, in the proofs of these theorems, we do not use that the new iterate $x_{k+1}$ is the exact solution of the proximal operator (see \Cref{proof:lemma:step2}). 
Instead, we use two properties of the proximal point operator:
\begin{itemize}
    \item decrease in function value:  $\phi_{k}(x_{k+1}) \leq \phi_{k}(x_{k})$;
    \item stationarity: $\nabla \phi_{k}(x_{k+1}) = 0$.
\end{itemize}
Thus, we can replace the computation of the exact proximal point in \Cref{alg:spam-exact} with finding a point that satisfies the above two conditions. 
Furthermore, we will relax the latter condition by taking an approximate stationary point. 
These arguments are summarized in the below assumption. 
\begin{definition}[\rm a\text{-}prox]\label{def:inexact}
    For a given client $k$, a gradient estimator $g_k$, a current state $x_k$, a stepsize $\gamma_k$ and a precision level $\epsilon$, the approximate proximal point $\aprox{\epsilon}{x_k,g_k,\gamma_k,k}$ is the set of vectors $y_{\sf ap}$, which satisfy
    % \vspace{-.3cm}
    \begin{itemize}
        \item decrease in function value:  $\Exp{\phi_k(y_{\sf ap}) } \leq \phi_k(x^k)$,%\vspace{-.3cm}
        \item approximate stationarity: $\Exp{\sqn{\nabla \phi_k(y_{\sf ap}) }} \leq \epsilon^2$,
    \end{itemize}
    % \vspace{-.3cm}
\end{definition}
where $\phi_k$ is defined in \eqref{eq:inexact}.
The pseudocode is described in \Cref{alg:spam-exact}.

\begin{theorem}[\algname{SPAM-inexact}]\label{thm:spam-inexact}
    Consider \algname{SPAM-inexact} for an objective function $f$ that satisfies Assumptions \ref{as:sigma} and \ref{as:similarity}. 
    Let $\gamma_k$ be a sequence of varying stepsizes satisfying $\gamma^2 \leq \frac{1}{16\delta^2}$ and choose $p_k = \frac{96\delta^2\gamma_k^2}{96\delta^2\gamma_k^2 + 1}$. 
    Then, 
    \begin{eqnarray}\label{eq:spam-decreasing-inexact-thm}
   \frac{1}{\Gamma_K}\sum_{k=1}^{K} \gamma_k\Exp{\sqn{\nabla f(x_{k+1})}} 
    &\leq &  \frac{40V_0}{\Gamma_K} 
    + \frac{2}{\Gamma_K}\sum_{k=1}^{K}{p_k\gamma_k^2}\sigma^2 
     + \frac{\epsilon^2}{8}.
    \end{eqnarray}
where $\Gamma_K = \sum_{k = 1}^{K} \gamma_k$. 
\end{theorem}
The proof is postponed to \Cref{app:proof-inexact}.
We observe that the level of inexactness $\epsilon^2$ appears explicitly in the theorem. In case when $\epsilon = 0$, we recover the result in \Cref{thm:spam-decaying} up to constants. 
\algname{SPAM-inexact} allows to avoid solving the local minimization problem required for finding the inexact proximal point operator. 
This is a significant improvement over \algname{SPAM}, as the latter requires minimizing (potentially) non-convex objectives at each iteration.

\section{Partial participation}\label{sec:spam-pp-new}

In this section, we present the most general form of our algorithm, which works with the approximate proximal operator and samples multiple clients (cohort) at each round. 
Specifically, it uses the random cohort $S_k$ to construct a better gradient estimator $g_k$. 
This gradient estimator is then broadcasted to a single random client $\xi_{k} \sim \cD$, who locally computes the approximate proximal point. 
The pseudocode can be found in \Cref{alg:spam-pp-new}.
\begin{algorithm}[h]
    \begin{algorithmic}[1]
        \STATE  \textbf{Input:}  learning rate $\gamma>0$, cohort size $B$, starting point $x_0\in\R^d$; \\proximal precision level $\epsilon$; initialize $g_0 = g_{-1}$;
        \FOR {$k=0,1,2, \ldots$}
        \STATE  \server samples a subset of clients $S_k$, with size $|S_k| = B$;
        \STATE  \server broadcasts $x_{k}$ to the clients from $S_k$.
        \FOR {$i \in S_k$ in parallel}
        \STATE \clients Set $g_k^{i} =  \nabla f_{i} (x_k) +(1-p_k)\left(g_{k-1} -\nabla f_{i} (x_{k-1}) \right)$;
        \STATE \clients Send $g_k^{i}$ to the server;
        \ENDFOR
        \STATE  \server $g_k = \frac{1}{B}\sum_{i \in S_k} g_k^{i}$ ;
        \STATE  \server {$\xi^{k+1} \sim \cD$};
        \STATE  \client $x_{k+1} \in \aprox{\epsilon}{x_k,g_k,\gamma_k,\xi^{k+1}}$;
        \ENDFOR
    \end{algorithmic}
    \caption{\algname{SPAM-PP}}
    \label{alg:spam-pp-new}    
\end{algorithm}

\begin{theorem}[\algname{SPAM-PP}]\label{thm:sppm-mvr-pp-new}

Suppose Assumptions \ref{as:sigma} and \ref{as:similarity} are satisfied. If $\xi_k \sim {\sf Unif}(S_k)$ at every iteration, then the iterates of \algname{SPAM-PP} with $\gamma_k \leq \frac{1}{4\delta}$ and $p_k = \frac{96\delta^2\gamma_k^2}{96\delta^2\gamma_k^2 + B^2}$ satisfy
\begin{eqnarray*}
\frac{1}{\Gamma_K}\sum_{k=0}^{K-1}\gamma_k\Exp{\sqn{\nabla f(x_{k+1})}}
&\leq& \frac{40}{\Gamma_K}\left(V_0 - \Exp{V_{K}}\right) + \frac{240}{\Gamma_K}\sum_{k=0}^{K-1}p_k\gamma_k \frac{\sigma^2}{B} + 7.5 \epsilon^2.
\end{eqnarray*}

\end{theorem}
The proof of the theorem is postponed to \Cref{app:proof-pp-new}.
When the client cohort size $B$ increases, the neighborhood shrinks. 
This is intuitive as when $B \rightarrow \infty$, we can access the exact objective $f$, and the neighborhood will vanish.

\begin{corollary}\label{cor:pp-complexity}
For properly chosen constant $\gamma_k = \gamma$ and momentum parameter $p_k = p$ 
% Choose constant stepsize $\gamma_k = \gamma$ and momentum parameter $p_k = p$ as
% \begin{equation*}
%     \gamma = \min \br{\frac{1}{\delta}, \br{\frac{BF}{2 \delta^ 2\sigma^2 K}}^{1/3}},
%     \qquad p = \max(\gamma^2\delta^2, 1/K).
% \end{equation*}
% Let us fix the parameters $\gamma_k = \gamma$ and $p_k = p$, such that $\gamma$ and $p$ satisfy the conditions of \Cref{thm:sppm-mvr-pp-new}. 
% Then, the 
communication complexity of \algname{SPAM-PP}, to obtain $\varepsilon$ error is of order 
\begin{eqnarray*}
    \mathcal{O} \br{\frac{\delta F}{\varepsilon} + \frac{\sigma^2}{B\varepsilon} + \frac{\delta \sigma F}{\sqrt{B}\varepsilon^{3/2}}}.
\end{eqnarray*}    

\end{corollary}

This result significantly improves upon the lower bound for $L$-smooth case
$\mathcal{O} \br{\frac{L F}{\varepsilon} + \frac{\sigma^2}{B\varepsilon} + \frac{L \sigma F}{\sqrt{B}\varepsilon^{3/2}}}$
    % \mathcal{O} \br{\frac{L + \sigma^2}{\varepsilon} + \frac{L \sigma}{\varepsilon^{3/2}}}
when $\delta \ll L$ \cite{arjevani2023lower}. 
We also observe that the complexity is an increasing function of $\delta$. Thus, our bound improves when data on different clients is similar. Moreover, increasing cohort size $B$ brings acceleration but to a certain level.
The proof of the corollary can be found in \Cref{sec:SPAM-PP_complexity}. 

In \Cref{sec:spam-pp}, we present another version of \algname{SPAM-PP}, called \algname{SPAM-PPA}, which uses the sampled cohort of clients to compute local proximal points. 
These points are then communicated to the server, and the new iterate is their average. 
Hence, the name \algname{SPAM-PP} with Averaging.

\section{Experiments}\label{sec:experiments}

To empirically validate our theoretical framework and its implications, we focus on a carefully controlled setting similar to \cite{khaled2022faster, lin2024stochastic}. Specifically, we consider a distributed ridge regression problem formulated in \eqref{eq:quad_problem}, which allows us to calculate and control the Hessian similarity \(\delta\). An essential advantage of this optimization problem is that the proximal operator has an explicit (closed-form) representation and can be computed precisely (up to machine accuracy). This allows us to isolate the effect of varying parameters on the method's performance. Appendix \ref{app:experiment} provides more details on the setup.
% We study the convergence properties of \algname{SPAM-inexact}.
% defined as
% Our main goal is to study the convergence properties of \algname{SPAM} algorithm.
% \begin{equation} \label{eq:quad_problem}
%     f(x) = \ExpSub{\xi}{\sqn{A_{\xi}x - y_{\xi}}} + \frac{\lambda}{2} \sqn{x},
% \end{equation}
% where $\xi$ is uniform random variable over $\{1, \dots, n\}$ for $n = 10, \lambda=0.1$. 
% We follow a similar to \cite{lin2024stochastic} procedure for synthetic data generation, which allows us to calculate and control Hessian similarity $\delta$.

% Namely, a random matrix $A_0 \in \mathbb{R}^{d\times d}$ ($d=100$) is generated with entries from a standard Gaussian distribution $\mathcal{N}(0, 1)$. Then we obtain $A = A_0 A_0^{\top}$ (to ensure symmetry) and set $A'_{\xi} = A + B_{\xi}$ by adding a random symmetric matrix $B_{\xi}$. Afterwards we modify $A_{\xi} = A'_{\xi} + I \lambda_{\min} (A'_{\xi})$ by adding an identity matrix $I$ times minimum eigenvalue to guarantee $A_{\xi} \succeq 0$. Entries of vectors $y_\xi \in \sR^d$, and initialization $x_0 \in \sR^d$ are generated from a standard Gaussian distribution $\mathcal{N}(0, 1)$.

% Our main goal in this section is to study the convergence properties of \algname{SPAM-inexact}. An important advantage of the considered problem formulation \eqref{eq:quad_problem} is that the proximal operator has an explicit (closed-form) representation and can be computed precisely (up to machine precision). This allows to isolate the effect of varying parameters on the method's performance. 

\begin{figure}[h]
\centering
\includegraphics[width=\linewidth]{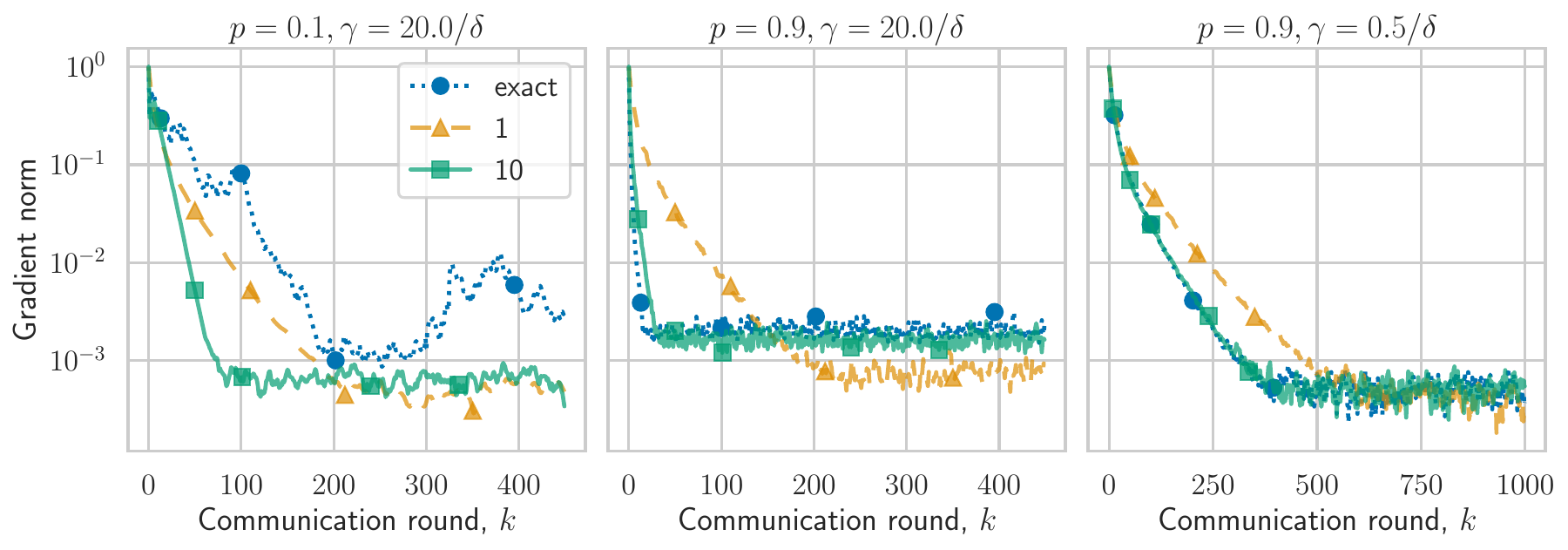}
\caption{Convergence of \algname{SPAM-inexact} on a ridge regression problem with different $p$ and $\gamma$.}
\label{fig:quadratic}
% \label{fig:varying_p}
\end{figure}

In Figure \ref{fig:quadratic}, we display convergence of Algorithm \ref{alg:spam-exact} with constant parameters $p$ and $\gamma$. The legend is shared, and labels refer to proximal operator computations: \say{exact} means using closed-form solution, \say{1} and \say{10} correspond to the number of local gradient descent steps. 
We evaluate the logarithm of a relative gradient norm $\log(\norm{\nabla f(x_k)}/\norm{\nabla f(x_0)})$ in the vertical axis. At every iteration, one client is sampled uniformly at random.

\textbf{Observations.} All the plots indicate convergence of the method to the neighborhood of the stationary point, followed by subsequent oscillations around the error floor.
The first ({left}) plot shows that for small momentum \( p = 0.1 \) and \(\gamma\) exceeding the theoretical bound \(1/\delta\), the algorithm can be very unstable with exact proximal point computations. Interestingly, approximate computation (1 or 10 local steps) results in more robust convergence.
The second ({middle}) plot demonstrates that a greater \( p = 0.9 \) results in steady convergence even for misspecified (too large) \(\gamma\). In addition, one can observe that in this case, more accurate proximal point evaluation results in significantly faster convergence but to a larger neighborhood than for one local step. This agrees well with observations for local gradient descent methods \cite{khaled2020tighter}.
The last ({right}) figure shows that a properly chosen, smaller \(\gamma = 0.5/\delta\) slows down convergence (twice as many communication rounds are shown). However, the method reaches a significantly lower error floor (as the vertical axis is shared across plots), which does not depend much on the accuracy of proximal point operator calculation. Moreover, 10 local steps are enough for basically the same fast convergence as with exact proximal point computation.

We would like to specifically note that momentum-based variance reduction has already shown empirical success \cite{karimireddy2021breaking, horvath2022fedshuffle} in practical federated learning scenarios. That is why our experiments focus on simpler but insightful setting to carefully study the properties of the proposed algorithm.

\section{Conclusion}\label{sec:conclusion}

We introduced \algname{SPAM}, an algorithm tailored for cross-device federated learning, which combines momentum variance reduction with the stochastic proximal point method.
Operating under second-order heterogeneity and bounded variance conditions, \algname{SPAM} does not necessitate smoothness of the objective function.
In its most general form, \algname{SPAM} achieves optimal communication complexity.
Furthermore, it does not prescribe a specific local method for analysis, providing practitioners with flexibility and responsibility in selecting suitable local solver.

\paragraph{Limitations and future work.}
The paper is of a theoretical nature and focuses on improving the understanding of stochastic non-convex optimization under Hessian similarity in the context of cross-device federated learning. We believe that separate experiments should be conducted to evaluate the experimental performance in a setting close to real life.

In standard optimization, the stepsize usually depends on the smoothness parameter. 
Adaptive methods allow iterative adjustment of the stepsize without additional information. 
In our case, the smoothness parameter is replaced by the second-order heterogeneity parameter $\delta$, on which the stepsize and momentum sequences of \algname{SPAM} depend. Removing this dependence using adaptive techniques under general assumptions remains an open problem even for the server-only \algname{MVR}, which serves as the basis for our algorithm.

Finally, federated learning comprises other aspects that we have not discussed above. These include privacy, security, personalization, etc., while our focus is on optimization and communication complexity. The study of their interplay is vital in future work. Our flexible framework can be beneficial due to the significantly simpler theory in contrast to prior works \cite{karimireddy2021breaking, patel2022towards}. % which can be more 
% We leave the study of their interplay as future work.

\section*{Acknowledgements}
We would like to thank Alexander Tyurin (KAUST) for useful discussions on Momentum Variance Reduction.

% \begin{figure}[t]
% \centering
% \includegraphics[width=\linewidth]{../plots/p_gamma_var.pdf}
% \caption{}
% \label{fig:varying_gamma}
% \end{figure}

% \begin{figure}[h]
% \centering
% \begin{subfigure}{0.5\textwidth}
%     \includegraphics[width=\linewidth]{plots/val_accuracy-a5a.pdf}
%     \caption{Validation accuracy of Gradient Descent ($\mathsf{GD}$) and Double Sketched Gradient Descent ($\mathsf{DSGD}$).}
%     \label{fig:val_accuracy}
% \end{subfigure}
% \begin{subfigure}{0.49\textwidth}
%     \includegraphics[width=\linewidth]{plots/test_acc_box-a5a.pdf}
%     \caption{Test accuracies distributions of sparsified solutions for the standard (ERM) formulation \eqref{eq:main} and MAST problem \eqref{eq:pretrained_compressed_problem}.}
%     \label{fig:test_acc_box}
% \end{subfigure}
% \caption{Validation and test performance of models based on ERM and MAST (with \rnd{$K$} sketches) approaches for different $q=K/d$.}
% \label{fig:val_test_acc}
% \end{figure}

% \clearpage

\bibliography{arxiv}
\bibliographystyle{alpha}

\clearpage
\addtocontents{toc}{\protect\setcounter{tocdepth}{2}}

%%%%%%%%%%%%%%%%%%%%%%%%%%%%%%%%%%%%%%%%%%%%%%%%%%%%%%%%%%%%

\appendix

\tableofcontents

\clearpage

\section{Convergence analysis for SPAM}\label{app:proof-spam}

\subsection{Proof of \Cref{prop:lyapunov}}\label{app:proof-lyapunov}
Recall that 
\begin{equation*}
    V_{k} = f(x_k) - f_{\inf} + \frac{3\gamma_k}{2p_k - p_k^2}\sqn{g_k - \nabla f(x_k)}.
\end{equation*}
We bound each term separately. 
We formulate three technical lemmas that are inspired by \cite{mishchenko2023federated}. 
Their proofs can be found in \Cref{app:proofs-lemmas}. 
We start by bounding the first term of the Lyapunov function, which is the function values.
\begin{lemma}\label{lemma:step2}
    Under the conditions of \Cref{prop:lyapunov}, the following recurrent inequality takes place
    \begin{equation}
\label{eq:step_2}
    f(x_{k+1}) - f_{\inf} \leq f(x_{k}) - f_{\inf} -\frac{1}{4\gamma_k}\sqn{x_{k+1}-x_k} + 2\gamma_k\sqn{\nabla f(x_k) - g_k}
\end{equation}
\end{lemma}

Then, we bound the second term of $V_k$. 
\begin{lemma}\label{lemma:step1}
    Under the conditions of \Cref{prop:lyapunov}, the following recurrent inequality takes place
    \begin{equation}\label{eq:step_1}
    \ExpCond{\sqn{g_{k+1} - \nabla f(x_{k+1})}}{\cF_k} 
    \leq (1-p_k)^2 \sqn{g_{k} -\nabla f(x_k)}  + 2(1-p_k)^2 \delta^2\sqn{ x_{k+1} - x_{k} } + 2p_k^2 \sigma^2.
    \end{equation}
\end{lemma}

 We observe that the term  $\sqn{x_{k+1}-x_k}$ is in both upper bounds. The following lemma provides a lower bound for this expression.

\begin{lemma}\label{lemma:step3}
Under the conditions of \Cref{prop:lyapunov}, the following recurrent inequality is true
 \begin{equation}
     \label{eq:step_3}
     \Exp{\sqn{x_{k+1}-x_k}} \geq \frac{\gamma_k^2}{4}\Exp{\sqn{\nabla f(x_{k+1})}} - \gamma_k^2\Exp{\sqn{g_k - \nabla f(x_k)}}. 
 \end{equation}
\end{lemma}

We now combine the results of the lemmas to bound $V_{K+1}$:
 \begin{eqnarray*}
     \Exp{V_{k+1}} &\overset{\eqref{eq:step_2} + \eqref{eq:step_1}}{\leq}& \alpha(1-p_k)^2 \sqn{g_{k} 
     - \nabla f(x_k)}  + 2\alpha\delta^2  (1-p_k)^2 \sqn{ x_{k+1} - x_{k} } + 2\alpha p_k^2 \sigma^2\\
     && + \Exp{f(x_{k}) - f_{\inf}} -\frac{1}{4\gamma_k}\Exp{\sqn{x_{k+1}-x_k}} + 2\gamma_k\Exp{\sqn{\nabla f(x_k) - g_k}}\\
     &=& \Exp{V_k} + \left(2\alpha\delta^2(1-p_k)^2-\frac{1}{4\gamma_k}\right)\Exp{\sqn{x_{k+1}-x_k}} + 2\alpha p_k^2\sigma^2 \\
     &&+ (2\gamma_k  - \alpha (2p_k - p_k^2))\Exp{\sqn{\nabla f(x_k) - g_k}}.
  \end{eqnarray*}
  The last inequality is true for every positive value of $\alpha$. 
  Let us now choose $\alpha = \frac{3\gamma_k}{2p_k - p_k^2}$. 
  Then,
  \begin{equation*}
       2\alpha\delta^2(1-p_k)^2-\frac{1}{4\gamma_k} 
       = \frac{6\gamma_k\delta^2(1-p_k)^2}{2p_k - p_k^2}-\frac{1}{4\gamma_k}  \leq - \frac{1}{8\gamma_k},
   \end{equation*} 
   where the latter is due to $4\delta\gamma_k \leq \sqrt{\nicefrac{p_k}{6(1 - p_k)}}$. 
   Therefore, we deduce 
  \begin{eqnarray*}
     \Exp{V_{k+1}}    
     % % &{=}& \Exp{V_k} + \left(3\delta^2\frac{1-p_k}{p}-\frac{1}{4\gamma_k}\right)\Exp{\sqn{x_{k+1}-x_k}} + 3\gamma_k\sigma^2 \\
     % && -\gamma_k\Exp{\sqn{\nabla f(x_k) - g_k}}\\
     &{\leq}& \Exp{V_k} - \frac{1}{8\gamma_k}\Exp{\sqn{x_{k+1}-x_k}} 
     - \gamma_k\Exp{\sqn{\nabla f(x_k) - g_k}} + 2\alpha p_k^2\sigma^2 \\
     &\overset{\eqref{eq:step_3}}{\leq}& \Exp{V_k} - \frac{\gamma_k}{32}\Exp{\sqn{\nabla f(x_{k+1})}} 
     + \frac{\gamma_k}{8}\Exp{\sqn{\nabla f(x_k) - g_k}} \\
     &&
     - \gamma_k\Exp{\sqn{\nabla f(x_k) - g_k}} + \frac{6\gamma_k p_k}{2 - p_k}\sigma^2\\
     &\leq &  \Exp{V_k} - \frac{\gamma_k}{32}\Exp{\sqn{\nabla f(x_{k+1})}} + 6\gamma_k p_k\sigma^2.
 \end{eqnarray*}
%  \esnote{I would write it a bit more explicitly for clarity.}
% where in $(\ast)$ we used 
% \begin{equation}
%     \gamma_k^2 \leq \frac{4p_k}{15\delta^2(1-p_k)} \leq \frac{2(2p_k - p_k^2)}{15\delta^2(1-p_k)^2} 
%     \Leftrightarrow \gamma_k \leq \frac{1}{8\delta^2(1-p_k)^2\alpha}.
% \end{equation}
This concludes the proof of the proposition.

\subsection{Proof of \Cref{thm:sppm-mvr}}\label{app:proof-thm:sppm-mvr}
Let us apply \Cref{prop:lyapunov} for the fixed stepsize $\gamma_k = \gamma$ and a fixed momentum coefficient $p_k = p$.
\begin{eqnarray*}
     \Exp{V_{k+1}}    
     &\leq &  \Exp{V_k} - \frac{\gamma}{32}\Exp{\sqn{\nabla f(x_{k+1})}} + 6\gamma  p\sigma^2.
 \end{eqnarray*}
Summing up these inequalities for $k = 0,\ldots, K-1$ leads to 
\begin{eqnarray*}
    \frac{1}{K}\sum_{k=1}^{K}\Exp{\sqn{\nabla f(x_{k})}} &\leq& \frac{32}{\gamma K}\left(V_0 - \Exp{V_{K}}\right) + 192 p\sigma^2\\
    &\leq& \frac{32(f(x_0) - f_{\inf})}{\gamma K} + \frac{30\sqn{g_0 - \nabla f(x_0)}}{ (2p - p^2)K} + 192 p\sigma^2.
\end{eqnarray*}
where $\gamma^2 \leq \min\left\{\frac{1}{16\delta^2},\frac{4p}{3\delta^2(1-p)}\right\}$.
This concludes the proof of the theorem.

\subsection{Proof of \Cref{thm:spam-decaying}} \label{app:proof-decaying}
From \Cref{prop:lyapunov} we have
\begin{eqnarray*}
     - \frac{\gamma_k}{32}\Exp{\sqn{\nabla f(x_{k+1})}}    
     &\leq &  \Exp{V_k} - \Exp{V_{k+1}}  + 6\gamma_k p_k\sigma^2.
 \end{eqnarray*}
Let us sum up these inequalities for $k = 0,1,\ldots,K-1$. 
We have a telescoping sum on the right-hand side. Then, dividing both sides on $\Gamma_K = \sum_{i=1}^{K} \gamma_i$, we deduce the following bound:
\begin{equation*}
    \frac{1}{\Gamma_K}\sum_{k=1}^{K}\gamma_k\Exp{\sqn{\nabla f(x_{k})}} 
    \leq \frac{32V_0}{\Gamma_K} 
    + \frac{2}{\Gamma_K}\sum_{k=1}^{K}\frac{15\delta^2\gamma_k^3}{15\delta^2\gamma_k^2 + 4} \sigma^2.
\end{equation*}
This concludes the proof.

\subsection{Proof of \Cref{thm:spam-inexact}}\label{app:proof-inexact}

We start by repeating the steps of the proof for \Cref{prop:lyapunov}. 
Notice that the proposition statement assumes that the iterate is exactly equal to the proximal point operator. 
However, as stated in \Cref{sec:inexact}, in the proofs of lemmas \ref{lemma:step2} and \ref{lemma:step1} we only use the property that $\phi_{k}(x_{k+1}) \leq \phi_{k}(x_{k})$ (see \eqref{eq:phi-decrease}).
Thus, both \eqref{eq:step_2} and \eqref{eq:step_1} are true for \algname{SPAM-inexact}. 
Therefore, 
\begin{eqnarray*}
     \Exp{V_{k+1}}    
     &{\leq}\quad \Exp{V_k}& - \frac{1}{8\gamma_k}\Exp{\sqn{x_{k+1}-x_k}} 
     - (2\gamma_k -\alpha (2p_k - p_k^2))\Exp{\sqn{\nabla f(x_k) - g_k}}\\  
     & & + 2\alpha p_k^2\sigma^2.
 \end{eqnarray*}
Below, reformulate the adaptation of \Cref{lemma:step3} for the inexact case to lower bound the second term on the right-hand side. 
\begin{lemma}\label{lemma:step3inexact} 
    Under the conditions of \Cref{prop:lyapunov}, we have the following bound 
     \begin{equation}\label{eq:step3_inexact}
     \Exp{\sqn{x_{k+1}-x_k}} \geq \frac{\gamma_k^2}{5}\Exp{\sqn{\nabla f(x_{k+1})}} - \gamma_k^2\Exp{\sqn{g_k - \nabla f(x_k)}} - \gamma_k^2\epsilon^2. 
 \end{equation}
\end{lemma}
The proof can be found in \Cref{proof:step3inexact}.  Thus,
\begin{eqnarray*}
    \Exp{V_{k+1}}    
    &\overset{\eqref{eq:step3_inexact}}{\leq}& \Exp{V_k} - \frac{\gamma_k}{40}\Exp{\sqn{\nabla f(x_{k+1})}} 
    + \frac{15\gamma_k p_k}{8(2 - p_k)}\sigma^2 + \frac{\gamma_k\epsilon^2}{8}\\
    && + \frac{\gamma_k}{8}\Exp{\sqn{\nabla f(x_k) - g_k}} 
    - (2\gamma_k -\alpha (2p_k - p_k^2))\Exp{\sqn{\nabla f(x_k) - g_k}}\\
    &\leq &  \Exp{V_k} - \frac{\gamma_k}{40}\Exp{\sqn{\nabla f(x_{k+1})}} + 2\gamma_k p_k\sigma^2 
    + \frac{\gamma_k\epsilon^2}{8}.
\end{eqnarray*}
Repeating this step for $k = 0,\ldots,K-1$, we deduce
\begin{eqnarray*}
   \frac{1}{\Gamma_K}\sum_{k=0}^{K-1} \gamma_k\Exp{\sqn{\nabla f(x_{k+1})}} 
    &\leq &  \frac{40V_0}{\Gamma_K} 
    + \frac{2}{\Gamma_K}\sum_{k=0}^{K-1} \frac{15\sigma^2\gamma_k^3}{15\sigma^2\gamma_k^2 + 4} \sigma^2 
     + \frac{\epsilon^2}{8}.
\end{eqnarray*}

\section{Proof of \Cref{thm:sppm-mvr-pp-new}}\label{app:proof-pp-new}

The proof follows the logic of \Cref{prop:lyapunov}. 
Recall that 
\begin{equation*}
    V_{k} = f(x_k) - f_{\inf} + \frac{3\gamma_k}{2p_k - p_k^2}\sqn{g_k - \nabla f(x_k)}.
\end{equation*}
Recall that \Cref{lemma:step2} is true for any gradient estimator $g_k$. 
Thus, \eqref{eq:step_2} is also valid for \algname{SPAM-PP}. 
Next, we estimate the second term of the Lyapunov function.
Recall that 
\begin{equation*}
\begin{aligned}
    g_{k+1} &= \frac{1}{S_k} \sum_{i \in S_k}  \curlybr{\nabla f_{i} (x_{k+1}) +(1-p_k)\left(g_{k} -\nabla f_{i} (x_{k}) \right)} \\
    & = \nabla \tilde{f}_k(x_{k+1}) + (1-p_k)\left(g_{k} -
    \nabla \tilde{f}_{k} (x_{k}) \right),
\end{aligned}
\end{equation*}
where  $\tilde{f}_k(x) := \frac{1}{S_k} \sum_{i\in S_k} \nabla f_{i}(x)$.  
Notice also that $\Exp{\tilde{f}_k(x)} = f(x)$, for every fixed $x\in\R^d$.
Furthermore, combining the convexity of the Euclidean norm and Hessian similarity \eqref{eq:similarity_almost_sure} we deduce that the estimator $\tilde{f}_k$ satisfies the Hessian similarity condition
\begin{align*}
    \norm{\nabla \tilde{f}_{k}(x) - \nabla f(x)  - \nabla \tilde{f}_{k}(y) + \nabla f(y)} 
    &\leq \frac{1}{B}\sum_{i\in S_k} \norm{\nabla {f}_{i}(x) - \nabla f(x)  - \nabla {f}_{i}(y) + \nabla f(y)} \\
    &\leq \frac{\delta}{B} \norm{x - y}.
\end{align*}
Finally, Jensen's inequality implies that $\tilde{f}_k$ satisfies the bounded variance condition as well:
\begin{equation*}
    \Exp{\norm{ \nabla \tilde{f}_{k}(x) -  \nabla f(x)}} \leq \sigma^2/B.
\end{equation*}
Repeating the analysis exactly as in the proof of \Cref{lemma:step1}, we obtain
\begin{eqnarray}\label{eq:step_1-pp-new}
    \Exp{\sqn{g_{k+1} - \nabla f(x_{k+1})}} 
    &\leq&  (1-p_k)^2\Exp{\sqn{g_{k} - \nabla f(x_{k})}}  \notag\\  &&+ 2(1-p_k)^2\frac{\delta^2}{B^2}\Exp{\sqn{x_{k+1}-x_k}} +\frac{2p_k^2\sigma^2}{B}. 
\end{eqnarray}
Let us now bound the Lyapunov function using \eqref{eq:step_2} and \eqref{eq:step_1-pp-new}:
 \begin{eqnarray*}
     \Exp{V_{k+1}} &{\leq}&\Exp{f(x_{k}) - f_{\inf}} 
    + {2\gamma_k}\Exp{\sqn{\nabla f(x_k) - g_k}}
       - \frac{1}{4\gamma_k} \Exp{ \normsq{ x_{k+1} - x_k}}\\
     && + \alpha(1-p_k)^2\Exp{\sqn{g_{k} - \nabla f(x_{k})}} + 2\alpha(1-p_k)^2\frac{\delta^2}{B^2}\Exp{\sqn{x_{k+1}-x_k}} +\frac{2\alpha p_k^2\sigma^2}{B} \\
     &=& \Exp{V_k} + \left(2\alpha\frac{\delta^2}{B^2}(1-p_k)^2-\frac{1}{4\gamma_k}\right)\Exp{\sqn{x_{k+1}-x_k}} + \frac{2\alpha p_k^2 \sigma^2}{B} \\
     &&+ (2\gamma_k - \alpha (2 p_k - p_k^2))\Exp{\sqn{\nabla f(x_k) - g_k}}.
 \end{eqnarray*}
 The latter is true for every positive $\alpha$. 
 Let us now plug in the value of $\alpha = \frac{3\gamma_k}{2p_k - p_k^2}$. 
 Then, using $\gamma \leq \sqrt{\frac{B^2p_k}{96\delta^2(1-p_k)}}$, we obtain
 \begin{equation}
    2\alpha\frac{\delta^2}{B^2}(1-p_k)^2-\frac{1}{4\gamma_k} 
    \leq    \frac{6\gamma_k\delta^2}{B^2(2p_k - p_k^2)}(1-p_k)^2-\frac{1}{4\gamma_k}
    \leq -\frac{1}{8\gamma_k}.
 \end{equation}
Hence, we have the following bound
 \begin{eqnarray*}
     \Exp{V_{k+1}} &{\leq}& 
     \Exp{V_k} -\frac{1}{8\gamma_k}\Exp{\sqn{x_{k+1}-x_k}} - \gamma_k \Exp{\sqn{\nabla f(x_k) - g_k}} + \frac{6 p_k\gamma_k \sigma^2}{B(2-p_k)} \\
     &\overset{\eqref{eq:step3_inexact}}{\leq}&\Exp{V_k} -\frac{1}{8\gamma_k}\br{\frac{\gamma_k^2}{5}\Exp{\sqn{\nabla f(x_{k+1})}} - \gamma_k^2\Exp{\sqn{g_k - \nabla f(x_k)}} - \gamma_k^2\epsilon^2 }\\ 
      &&- \gamma_k \Exp{\sqn{\nabla f(x_k) - g_k}} + \frac{6 p_k\gamma_k \sigma^2}{B} \\
     &{\leq}& \Exp{V_k} - \frac{\gamma_k}{40}\Exp{\sqn{\nabla f(x_{k+1})}} - \frac{7\gamma_k}{8} \Exp{\sqn{\nabla f(x_k) - g_k}} + \frac{6 p_k\gamma_k \sigma^2}{B} + \frac{\gamma_k \epsilon^2}{8}\\
     &\leq&  \Exp{V_k} - \frac{\gamma_k}{40}\Exp{\sqn{\nabla f(x_{k+1})}}  + \frac{6 p_k\gamma_k \sigma^2}{B} + \frac{\gamma_k \epsilon^2}{8}.
 \end{eqnarray*}

Thus, we have 
\begin{eqnarray*}
    \frac{1}{\Gamma_K}\sum_{k=0}^{K-1}\gamma_k\Exp{\sqn{\nabla f(x_{k+1})}} 
    &\leq& \frac{40}{\Gamma_K}\left(V_0 - \Exp{V_{K}}\right) + \frac{240}{\Gamma_K}\sum_{k=0}^{K-1}p_k\gamma_k \frac{\sigma^2}{B} + {7.5 \epsilon^2}.
\end{eqnarray*}
This concludes the proof of the theorem.

\section{Complexity analysis of the methods} 
We use $\lesssim$ to ignore numerical constants in the subsequent analysis. 

\subsection{Proof of \Cref{cor:spam-complexity}}\label{app:cor:spam-complexity}
We have stepsize condition $\gamma \lesssim \min\{1/\delta, \sqrt{p/(\delta^2(1-p))}\}$, which implies that $\gamma \lesssim\sqrt{p}/\delta$ or $p \gtrsim (\gamma \delta)^2$.
Denote $F \eqdef f(x_0) - f_{\inf}$, then convergence rate of \algname{SPAM} can be expressed as
% takes the following form.
% \begin{equation}
%     \frac{1}{K}\sum_{k=1}^{K} \Exp{\sqn{\nabla f(x_{k})}} \lesssim
%     \frac{\Delta}{\gamma K} + \frac{\sqn{g_0 - \nabla f(x_0)}}{ (2p - p^2)K} + p\sigma^2.
% \end{equation}
\begin{eqnarray*}
    R_K \eqdef \frac{1}{K}\sum_{k=1}^{K}\Exp{\sqn{\nabla f(x_{k})}} &\lesssim& 
    \frac{f(x_0) - f_{\inf}}{\gamma K} + \frac{\sqn{g_0 - \nabla f(x_0)}}{ (2p - p^2)K} + p\sigma^2 \\
    &\lesssim& 
    %\frac{\delta F}{K} + \frac{\delta F}{\sqrt{p} K} 
    \frac{F}{\gamma K} + \frac{\sqn{g_0 - \nabla f(x_0)}}{p K} + p\sigma^2,
\end{eqnarray*}
where in the last inequality, we used condition for the stepsize and the fact that $p (2-p) \geq p$.
Next, by using an argument similar to that in \cite{karimireddy2021breaking}, we suppose (without loss of generality) that the method is run for $K$ iterations. For the first $K/2$ iterations, we simply sample $\nabla f_{\xi}$ at $x_0$ to set $g_0 = \frac{1}{K/2} \sum_{i=1}^{K/2} \nabla f_{\xi_i}(x_0)$.
Then, according to \eqref{eq:sigma},  we have $\Exp{\sqn{g_0 - \nabla f(x_0)}} \leq \frac{\sigma^2}{K/2}.$
Now, choose $p = \max(\gamma^2\delta^2, 1/K)$
\begin{eqnarray*}
    R_K \lesssim
    % \frac{\delta F}{K} + \frac{\delta F}{\sqrt{p} K} 
    \frac{F}{\gamma K} + \frac{\sigma^2}{p K^2} + p\sigma^2 
    \lesssim
    \frac{F}{\gamma K} + \frac{\sigma^2}{K} + \gamma^2\delta^2\sigma^2 + \frac{\sigma^2}{K}.
\end{eqnarray*}
Next set $\gamma = \min \br{\frac{1}{\delta}, \br{\frac{F}{2 \delta^ 2\sigma^2 K}}^{1/3}}$ and the rate results in
\begin{eqnarray*}
    R_K &\lesssim&
    \frac{\delta F}{K} + \frac{F}{K}\br{\frac{2 \delta^ 2\sigma^2 K}{F}}^{1/3} + \br{\frac{F}{2 \delta^ 2\sigma^2 K}}^{2/3} \delta^2 \sigma^2 + \frac{\sigma^2}{K}
    \\ &\lesssim&
    \frac{\delta F + \sigma^2}{K} + \br{\frac{\delta \sigma F}{K}}^{2/3},
\end{eqnarray*}
which leads to the communication complexity of
$\mathcal{O} \br{\frac{\delta F + \sigma^2}{\varepsilon} + \frac{\delta \sigma F}{\varepsilon^{3/2}}}$.
This concludes the proof.

\subsection{Proof of \Cref{cor:pp-complexity}} \label{sec:SPAM-PP_complexity}
In this part, we analyze communication complexity similar to Section \ref{sec:SPAM-PP_complexity}. The focus is on constant stepsize case $\gamma_k \equiv \gamma \lesssim \min\{1/\delta, \sqrt{p} B/\delta\}$ and exact proximal computation $\epsilon=0$.
% Note that we focus on communication complexity which is typically the main bottleneck in cross-device federated setting \cite{kairouz2021advances}. We use $\lesssim$ to ignore numerical constants in the subsequent analysis. stepsize condition $\gamma \lesssim \min\{1/\delta, \sqrt{p/(\delta^2(1-p))}\}$, which implies that $\gamma \lesssim\sqrt{p}/\delta$ or $p \gtrsim (\gamma \delta)^2$.
Denote $F \eqdef f(x_0) - f_{\inf}$, then convergence rate of \algname{SPAM-PP} can be expressed as
\begin{eqnarray*}
    R_K \eqdef \frac{1}{K}\sum_{k=1}^{K}\Exp{\sqn{\nabla f(x_{k})}} &\lesssim& 
    \frac{F}{\gamma K} + \frac{\sqn{g_0 - \nabla f(x_0)}}{p K} + \frac{p\sigma^2}{B}.
\end{eqnarray*}

By using the same reasoning as in \ref{app:cor:spam-complexity} set 
$$g_0 = \frac{1}{B K/2} \sum_{j=1}^{K/2} \sum_{i=1}^{S_j} \nabla f_{\xi_i}(x_0)$$ to make sure
$\Exp{\sqn{g_0 - \nabla f(x_0)}} \leq \sigma^2/(B K/2)$. Then
\begin{eqnarray*}
    R_K \lesssim \frac{F}{\gamma K} + \frac{\sigma^2}{p BK^2} + \frac{p\sigma^2}{B}.
\end{eqnarray*}

Now choose $p = \max(\gamma^2\delta^2, 1/K)$ which leads to
\begin{eqnarray*}
    R_K \lesssim \frac{F}{\gamma K} + \frac{\sigma^2}{BK} + \frac{\gamma^2\delta^2\sigma^2}{B} + \frac{\sigma^2}{BK}.
\end{eqnarray*}
Next set $\gamma = \min \br{\frac{1}{\delta}, \br{\frac{BF}{2 \delta^ 2\sigma^2 K}}^{1/3}}$ and the rate results in
\begin{eqnarray*}
    R_K &\lesssim&
    \frac{\delta F}{K} + \frac{F}{K}\br{\frac{2 \delta^ 2\sigma^2 K}{BF}}^{1/3} + \frac{\sigma^2}{BK} + \br{\frac{BF}{2 \delta^ 2\sigma^2 K}}^{2/3} \frac{\delta^2 \sigma^2}{B} + \frac{\sigma^2}{BK}
    \\ &\lesssim&
    \frac{\delta F}{K} + \frac{\sigma^2}{BK} + \br{\frac{\delta \sigma F}{\sqrt{B}K}}^{2/3},
\end{eqnarray*}
which leads to communication complexity 
\begin{eqnarray*}
    \mathcal{O} \br{\frac{\delta F}{\varepsilon} + \frac{\sigma^2}{B\varepsilon} + \frac{\delta \sigma F}{\sqrt{B}\varepsilon^{3/2}}}.
\end{eqnarray*}

\section{Partial participation with averaging}\label{sec:spam-pp}

\begin{algorithm}[h]
    \begin{algorithmic}[1]
        \STATE  \textbf{Input:}  learning rate $\gamma>0$, starting point $x_0\in\R^d$; \\proximal precision level $\epsilon$; initialize $g_0 = g_{-1}$;
        \FOR {$k=0,1,2, \ldots$}
        \STATE Sample a subset of clients $S_k$, with size $|S_k| = B$;
        \STATE Selected clients do local \algname{SPAM} updates;
        \FOR {$i \in S_k$}
        \STATE Set $g_k^{i} =  \nabla f_{i} (x_k) +(1-p_k)\left(g_{k-1} -\nabla f_{i} (x_{k-1}) \right)$;
        \STATE Broadcast $g_k^{i}$ to the server;
        \ENDFOR
        \STATE  $g_k = \frac{1}{B}\sum_{i \in S_k} g_k^{i}$ ;
        \FOR {$i \in S_k$}
        \STATE Set $x_{k+1}^{i} = \aprox{\epsilon}{x_k,g_k,\gamma_k,i}$;
        \STATE Broadcast $x_{k+1}^{i}$ to the server;
        \ENDFOR
        \STATE The server aggregates the local iterates: 
        $x_{k+1} = \frac{1}{B}\sum_{i\in S_k}x_{k+1}^{i}$;
        \ENDFOR
    \end{algorithmic}
    \caption{\algname{SPAM-PPA}}
    \label{alg:spam-pp}    
\end{algorithm}

\begin{theorem}[\algname{SPAM-PPA}]\label{thm:sppm-mvr-pp}
Suppose Assumptions  \ref{as:sigma}, \ref{as:similarity} are satisfied and the objective function $f$ is $L$-smooth. If $\xi_k \sim {\sf Unif}(S_k)$ at every iteration, then the iterates of \algname{SPAM-PPA} with $\gamma_k \leq \frac{1}{4(\delta + L)}$ and $p_k = \frac{96\delta^2\gamma_k^2}{96\delta^2\gamma_k^2 + B^2}$ satisfy
\begin{eqnarray*}
    \frac{1}{\Gamma_K}\sum_{k=0}^{K-1}\gamma_k\Exp{\sqn{\nabla f(x^{\xi}_{k+1})}} 
    &\leq& \frac{40}{\Gamma_K}\left(V_0 - \Exp{V_{K}}\right) + \frac{240}{\Gamma_K}\sum_{k=0}^{K-1}p_k\gamma_k \frac{\sigma^2}{B} + {7.5 \epsilon^2}.
\end{eqnarray*}
\end{theorem}
    
    The result of the theorem is similar to the one in \Cref{thm:sppm-mvr-pp-new}. 
    In fact, following the proof scheme of \Cref{cor:pp-complexity}, one can derive the complexity analysis for \algname{SPAM-PPA}.
    However, unlike previous results, we require the objective function $f$ to be smooth. 
    % This is a major drawback, as in practice, we do not have smoothness.

\subsection{Proof of \Cref{thm:sppm-mvr-pp}}\label{app:proof-pp}

The proof follows the logic of \Cref{prop:lyapunov}. 
Recall that 
\begin{equation*}
    V_{k} = f(x_k) - f_{\inf} + \frac{3\gamma_k}{2p_k - p_k^2}\sqn{g_k - \nabla f(x_k)}.
\end{equation*}
We start with proving a descent lemma. Recall  that $\xi_k \sim {\sf Unif}(S_k)$, for the fixed $S_k$.
\begin{lemma}\label{lemma:step_2-pp}
    For an $L$-smooth objective $f$ satisfying assumptions \ref{as:sigma},\ref{as:similarity} and parameters $\gamma_k^2 \leq \min\left\{\frac{1}{16(L + \delta)^2},\frac{4p_k}{15\delta^2(1-p_k)}\right\}$, the iterates of the  \algname{SPAM-PPA} algorithm satisfy 
    \begin{equation}
    \label{eq:step_2-pp}
        \Exp{f(x_{k+1}) - f_{\inf}} \leq \Exp{f(x_{k}) - f_{\inf}} 
        + {2\gamma_k}\Exp{\sqn{\nabla f(x_k) - g_k}}
           - \frac{1}{4\gamma_k} \Exp{ \normsq{ x^{\xi_k}_{k+1} - x_k}}.
    \end{equation}
\end{lemma}
The proof of the lemma is deferred to \Cref{app:proof-lemma:step_2-pp}.
Next, we estimate the second term of the Lyapunov function.
Recall that 
\begin{equation*}
\begin{aligned}
    g_{k+1} &= \frac{1}{S_k} \sum_{i \in S_k}  \curlybr{\nabla f_{i} (x_{k+1}) +(1-p_k)\left(g_{k} -\nabla f_{i} (x_{k}) \right)} \\
    & = \nabla \tilde{f}_k(x_{k+1}) + (1-p_k)\left(g_{k} -
    \nabla \tilde{f}_{k} (x_{k}) \right),
\end{aligned}
\end{equation*}
where  $\tilde{f}_k(x) := \frac{1}{S_k} \sum_{i\in S_k} \nabla f_{i}(x)$.  
Notice that $\Exp{\tilde{f}_k(x)} = f(x)$, for every fixed $x\in\R^d$.
Furthermore, combining the convexity of the Euclidean norm and Hessian similarity \eqref{eq:similarity_almost_sure} we deduce that the estimator $\tilde{f}_k$ satisfies the Hessian similarity condition
\begin{align*}
    \norm{\nabla \tilde{f}_{k}(x) - \nabla f(x)  - \nabla \tilde{f}_{k}(y) + \nabla f(y)} 
    &\leq \frac{1}{B}\sum_{i\in S_k} \norm{\nabla {f}_{i}(x) - \nabla f(x)  - \nabla {f}_{i}(y) + \nabla f(y)} \\
    &\leq \frac{\delta}{B} \norm{x - y}.
\end{align*}
Furthermore, Jensen's inequality implies that $\tilde{f}_k$ satisfies the bounded variance condition as well:
\begin{equation*}
    \Exp{\norm{ \nabla \tilde{f}_{k}(x) -  \nabla f(x)}} \leq \sigma^2/B.
\end{equation*}
Repeating the analysis exactly as in the proof of \Cref{lemma:step1}, we obtain
\begin{eqnarray*}
    \Exp{\sqn{g_{k+1} - \nabla f(x_{k+1})}} 
    &\leq&  (1-p_k)^2\Exp{\sqn{g_{k} - \nabla f(x_{k})}} \\ 
    &&+ 2(1-p_k)^2\frac{\delta^2}{B^2}\Exp{\sqn{x_{k+1}-x_k}} + \frac{2p_k^2\sigma^2}{B}.  
\end{eqnarray*}
Assume now that $\xi_k \sim {\sf Unif}(S_k)$, for a fixed $S_k$. 
The latter means $x_{k+1} = \ExpCond{x^{\xi_k}_{k+1}}{\cG_k}$, and subsequently, Jensen's inequality yields
\begin{equation}\label{eq:step_1-pp}
\begin{aligned}
    \Exp{\sqn{g_{k+1} - \nabla f(x_{k+1})}} 
    {\leq}& (1-p_k)^2\Exp{\sqn{g_{k} - \nabla f(x_{k})}} \\ 
    &+ 2(1-p_k)^2\frac{\delta^2}{B^2}\Exp{\sqn{\ExpCond{x^{\xi_k}_{k+1}}{\cG_k}-x_k}} +\frac{2p_k^2 \sigma^2}{B} \notag \\
    {\leq}& (1-p_k)^2\Exp{\sqn{g_{k} - \nabla f(x_{k})}} \\ 
    &+ 2(1-p_k)^2\frac{\delta^2}{B^2}\Exp{\ExpCond{\sqn{x^{\xi_k}_{k+1}-x_k}}{\cG_k}} +\frac{2p_k^2\sigma^2}{B} \notag \\
    =& (1-p_k)^2\Exp{\sqn{g_{k} - \nabla f(x_{k})}} \\ 
    &+ 2(1-p_k)^2\frac{\delta^2}{B^2}\Exp{\sqn{x^{\xi_k}_{k+1}-x_k}} +\frac{2p_k^2\sigma^2}{B}.
\end{aligned}
\end{equation}

Now, we need to bound $\Exp{ \normsq{ x^{\xi_k}_{k+1} - x_k}}$ from below. 

\begin{lemma}\label{lemma:step_3-pp}
    Under assumptions \ref{as:sigma} and \ref{as:similarity}, we have the following lower bound for the iterates of \algname{SPAM-PPA} algorithm
    \begin{equation}
     \label{eq:step_3-pp}
         \Exp{\sqn{x^{\xi_k}_{k+1}-x_k}} \geq \frac{\gamma_k^2}{5}\Exp{\sqn{\nabla f(x^{\xi_k}_{k+1})}} - \gamma_k^2\Exp{\sqn{g_k - \nabla f(x_k)}} - \gamma_k\epsilon^2. 
     \end{equation}
\end{lemma}
The proof of the lemma can be found in \Cref{app:proof-lemma:step_3-pp}.
Let us now bound the Lyapunov function using \eqref{eq:step_2-pp} and \eqref{eq:step_1-pp}:
 \begin{eqnarray*}
     \Exp{V_{k+1}} &{\leq}&\Exp{f(x_{k}) - f_{\inf}} 
    + {2\gamma_k}\Exp{\sqn{\nabla f(x_k) - g_k}}
       - \frac{1}{4\gamma_k} \Exp{ \normsq{ x^{\xi_k}_{k+1} - x_k}}\\
     && + \alpha(1-p_k)^2\Exp{\sqn{g_{k} - \nabla f(x_{k})}} + 2\alpha(1-p_k)^2\frac{\delta^2}{B^2}\Exp{\sqn{x^{\xi_k}_{k+1}-x_k}} +\frac{2\alpha p_k^2\sigma^2}{B} \\
     &=& \Exp{V_k} + \left(2\alpha\frac{\delta^2}{B^2}(1-p_k)^2-\frac{1}{4\gamma_k}\right)\Exp{\sqn{x^{\xi_k}_{k+1}-x_k}} + \frac{2\alpha p_k^2 \sigma^2}{B} \\
     &&+ (2\gamma_k - \alpha (2 p_k - p_k^2))\Exp{\sqn{\nabla f(x_k) - g_k}}.
 \end{eqnarray*}
 The latter is true for every positive $\alpha$. 
 Let us now plug in the value of $\alpha = \frac{3\gamma_k}{2p_k - p_k^2}$. 
 Then, using $\gamma \leq \sqrt{\frac{B^2p_k}{96\delta^2(1-p_k)}}$, we obtain
 \begin{equation}
    2\alpha\frac{\delta^2}{B^2}(1-p_k)^2-\frac{1}{4\gamma_k} 
    \leq    \frac{6\gamma_k\delta^2}{B^2(2p_k - p_k^2)}(1-p_k)^2-\frac{1}{4\gamma_k}
    \leq -\frac{1}{8\gamma_k}.
 \end{equation}
Hence, we have the following bound
 \begin{eqnarray*}
     \Exp{V_{k+1}} &{\leq}& 
     \Exp{V_k} -\frac{1}{8\gamma_k}\Exp{\sqn{x^{\xi_k}_{k+1}-x_k}} - \gamma_k \Exp{\sqn{\nabla f(x_k) - g_k}} + \frac{6 p_k\gamma_k \sigma^2}{B(2-p_k)} \\
     &\overset{\eqref{eq:step_3-pp}}{\leq}&\Exp{V_k} -\frac{1}{8\gamma_k}\br{\frac{\gamma_k^2}{5}\Exp{\sqn{\nabla f(x^{\xi_k}_{k+1})}} - \gamma_k^2\Exp{\sqn{g_k - \nabla f(x_k)}} - \gamma_k^2\epsilon^2 }\\ 
      &&- \gamma_k \Exp{\sqn{\nabla f(x_k) - g_k}} + \frac{6 p_k\gamma_k \sigma^2}{B} \\
     &{\leq}& \Exp{V_k} - \frac{\gamma_k}{40}\Exp{\sqn{\nabla f(x^{\xi_k}_{k+1})}} - \frac{7\gamma_k}{8} \Exp{\sqn{\nabla f(x_k) - g_k}} + \frac{6 p_k\gamma_k \sigma^2}{B} + \frac{\gamma_k \epsilon^2}{8}\\
     &\leq&  \Exp{V_k} - \frac{\gamma_k}{40}\Exp{\sqn{\nabla f(x^{\xi_k}_{k+1})}}  + \frac{6 p_k\gamma_k \sigma^2}{B} + \frac{\gamma_k \epsilon^2}{8}.
 \end{eqnarray*}

Thus, we have 
\begin{eqnarray*}
    \frac{1}{\Gamma_K}\sum_{k=0}^{K-1}\gamma_k\Exp{\sqn{\nabla f(x^{\xi_k}_{k+1})}} 
    &\leq& \frac{40}{\Gamma_K}\left(V_0 - \Exp{V_{K}}\right) + \frac{240}{\Gamma_K}\sum_{k=0}^{K-1}p_k\gamma_k \frac{\sigma^2}{B} + {7.5 \epsilon^2}.
\end{eqnarray*}
This concludes the proof of the theorem.

\section{Proofs of the technical lemmas}\label{app:proofs-lemmas}

%%%%%%%%%%%%%%%%%%%%%%%%%%%%%%%%%%%%%%%%%%%%%%%%%%%%%%%%%%%%
\subsection{Proof of \Cref{lemma:step2}}\label{proof:lemma:step2}

By the main theorem of Calculus, we have 
\begin{eqnarray*}
    f(x_{k+1}) - f(x_k) &=& \int^1_0\abr{\nabla f(\underbrace
    {x_k +\tau(x_{k+1}-x_k)}_{\eqdef x(\tau)}), x_{k+1} - x_k}d\tau,\\
    f_{\xi_k}(x_{k+1}) - f_{\xi_k}(x_k) &=& \int^1_0\abr{\nabla f_{\xi_k}(\underbrace
    {x_k +\tau(x_{k+1}-x_k)}_{\eqdef x(\tau)}), x_{k+1} - x_k}d\tau
\end{eqnarray*}
Therefore the difference in function value can be bounded as follows:
\begin{eqnarray*}
    f(x_{k+1}) - f(x_k) &=& f_{\xi_k}(x_{k+1}) - f_{\xi_k}(x_k) \\
    &&+\int^1_0\abr{\nabla f(x(\tau)) -\nabla f_{\xi_k}(x(\tau)) , x_{k+1} - x_k}d\tau\\
    &=& f_{\xi_k}(x_{k+1}) - f_{\xi_k}(x_k) + \abr{g_k - \nabla f_{\xi_k}(x_k), x_{k+1}-x_k}\\
    &&+\int^1_0\abr{\nabla f(x(\tau)) -\nabla f_{\xi_k}(x(\tau)) - g_k + \nabla f_{\xi_k}(x_k) , x_{k+1} - x_k}d\tau\\
    &\leq& -\frac{1}{2\gamma_k}\sqn{x_{k+1}-x_k} + \abr{\nabla f(x_k) - g_k, x_{k+1} - x_k} \\
    && + \int^1_0\abr{\nabla f(x(\tau)) -\nabla f_{\xi_k}(x(\tau)) - \nabla f(x_k) + \nabla f_{\xi_k}(x_k) , x_{k+1} - x_k}d\tau.
\end{eqnarray*}

The last inequality is due to
\begin{equation}\label{eq:phi-decrease}
    f_{\xi_k}(x_{k+1})+\abr{g_k-\nabla f_{\xi_k} (x_k), x_{k+1} - x_k}+ \frac{1}{2\gamma_k}\|x_{k+1}-x_k\|^2 \leq f_{\xi_k}(x_{k}), 
\end{equation}
which is a direct consequence of $x_{k+1} = \arg\min\limits_{x}\left\{f_{\xi_k}(x)+\abr{ g_k-\nabla f_{\xi_k} (x_k), x - x_k} + \frac{1}{2\gamma_k}\|x-x_k\|^2\right\}$. 
Let us now apply Cauchy-Schwartz inequality to bound both scalar products:
\begin{eqnarray*}
    f(x_{k+1}) - f(x_k)
    &\leq& -\frac{1}{2\gamma_k}\sqn{x_{k+1}-x_k} + 2\gamma_k\sqn{\nabla f(x_k) - g_k} + \frac{1}{8\gamma_k}\sqn{x_{k+1} - x_k}\\
    && + \int^1_0\norm{\nabla f(x(\tau)) -\nabla f_{\xi_k}(x(\tau)) - \nabla f(x_k) + \nabla f_{\xi_k}(x_k) }\norm{ x_{k+1} - x_k}d\tau\\
    &\overset{\eqref{eq:similarity_almost_sure}}{\leq}& -\frac{1}{2\gamma_k}\sqn{x_{k+1}-x_k} + 2\gamma_k\sqn{\nabla f(x_k) - g_k} + \frac{1}{8\gamma_k}\sqn{x_{k+1} - x_k}\\
    && + \delta \int^1_0\norm{x(\tau) - x_k }\norm{ x_{k+1} - x_k}d\tau\\
    &=& -\frac{1}{2\gamma_k}\sqn{x_{k+1}-x_k} + 2\gamma_k\sqn{\nabla f(x_k) - g_k} + \frac{1}{8\gamma_k}\sqn{x_{k+1} - x_k}\\
    && + \frac{\delta}{2} \sqn{ x_{k+1} - x_k}\\
    &\overset{\gamma_k \leq \frac{1}{4\delta}}{\leq} & -\frac{1}{4\gamma_k}\sqn{x_{k+1}-x_k} + 2\gamma_k\sqn{\nabla f(x_k) - g_k}.
\end{eqnarray*}
Thus, we have 
\begin{equation}
    f(x_{k+1}) - f_{\inf} \leq f(x_{k}) - f_{\inf} -\frac{1}{4\gamma_k}\sqn{x_{k+1}-x_k} + 2\gamma_k\sqn{\nabla f(x_k) - g_k}.
\end{equation}

This concludes the proof of the lemma.

\subsection{Proof of \Cref{lemma:step1}}

Recall that $g_{k+1} =  \nabla f_{\xi_{k+1}} (x_{k+1}) +(1-p_k)\left(g_{k} -\nabla f_{\xi_{k+1}} (x_{k}) \right)$. We define $\cF_k:= \{x_{k+1}, x_{k}, g_k\}$. Then, 
\begin{eqnarray*}
    &&\ExpCond{\sqn{g_{k+1} - \nabla f(x_{k+1})}}{\cF_k} \\
    \quad &=& \ExpCond{\sqn{\nabla f_{\xi_{k+1}} (x_{k+1}) +(1-p_k)\left(g_{k} -\nabla f_{\xi_{k+1}} (x_{k}) \right) - \nabla f(x_{k+1})}}{\cF_k}\\
    \quad &=& 
    \ExpCond{\sqn{\nabla f_{\xi_{k+1}} (x_{k+1}) - \nabla f (x_{k+1}) 
    + (1-p_k)\left(\nabla f(x_k)  -\nabla f_{\xi_{k+1}} (x_{k}) \right) }}{\cF_k} \\
    &&+ (1-p_k)^2 \sqn{g_{k} -\nabla f(x_k)}.
\end{eqnarray*}
The last equality is due to the bias-variance formula and the fact that $\xi_{k+1}$ is independent of $\cF_k$ and that the stochastic gradients are unbiased. 
Using the Cauchy-Schwartz inequality, we deduce the following bound for the first term on the right-hand side, where $\alpha > 0$ is an arbitrary constant:
\begin{align}\label{eq:step1-sigma2}
   \mathrm{E}&\sbr{\sqn{\nabla f_{\xi_{k+1}} (x_{k+1}) - \nabla f (x_{k+1}) 
    + (1-p_k)\left(\nabla f(x_k)  -\nabla f_{\xi_{k+1}} (x_{k}) \right) } \mid {\cF_k}} \notag \\
     & \quad =\mathrm{E}\left[\|p_k\br{\nabla f_{\xi_{k+1}} (x_{k+1}) - \nabla f (x_{k+1})}\right. \\
      & \qquad +  \left. (1-p_k)\left(\nabla f_{\xi_{k+1}} (x_{k+1}) - \nabla f (x_{k+1}) + \nabla f(x_k)  -\nabla f_{\xi_{k+1}} (x_{k}) \right) \|^2 \mid {\cF_k} \right] \notag \\
     & \quad \leq   (1 + \alpha)p_k^2 
    \ExpCond{\sqn{\nabla f_{\xi_{k+1}} (x_{k+1}) - \nabla f (x_{k+1})}}{\cF_k} \\
      & \qquad +  (1 + \alpha^{-1})(1-p_k)^2 \ExpCond{\sqn{\nabla f_{\xi_{k+1}} (x_{k+1}) - \nabla f (x_{k+1}) + \nabla f(x_k)  -\nabla f_{\xi_{k+1}} (x_{k})  }}{\cF_k}. \notag
\end{align}
We apply \eqref{eq:sigma} and \eqref{eq:similarity_almost_sure} to bound, respectively, the first term and the second term on the right-hand side of \eqref{eq:step1-sigma2}:
\begin{eqnarray*}
    &&\ExpCond{\sqn{\nabla f_{\xi_{k+1}} (x_{k+1}) - \nabla f (x_{k+1}) 
    + (1-p_k)\left(\nabla f(x_k)  -\nabla f_{\xi_{k+1}} (x_{k}) \right) }}{\cF_k} \\
    & \qquad \leq&  (1 + \alpha)p_k^2 \sigma^2 + (1 + \alpha^{-1})(1-p_k)^2 \delta^2\sqn{ x_{k+1} - x_{k} }. 
\end{eqnarray*}
Taking $\alpha = 1$, we obtain the following
\begin{equation*}
    \ExpCond{\sqn{g_{k+1} - \nabla f(x_{k+1})}}{\cF_k} 
    \leq (1-p_k)^2 \sqn{g_{k} -\nabla f(x_k)}  + 2(1-p_k)^2 \delta^2\sqn{ x_{k+1} - x_{k} } + 2p_k^2 \sigma^2.
\end{equation*}
This concludes the proof of the lemma.
\subsection{Proof of \Cref{lemma:step3}}

By the definition of $x_{k+1}$, we have 
\begin{eqnarray*}
    \sqn{x_{k+1}-x_k} &=& \gamma_k^2 \sqn{\nabla f_{\xi_k}(x_{k+1}) + g_k - \nabla f_{\xi_k}(x_{k}) }\\
    &=& \gamma_k^2 \sqn{ \nabla f(x_{k+1}) +g_k - \nabla f(x_k) + \nabla f_{\xi_k}(x_{k+1}) -\nabla f(x_{k+1}) - \nabla f_{\xi_k}(x_{k}) + \nabla f(x_k) }\\
    &\geq&\frac{\gamma_k^2}{3}\sqn{\nabla f(x_{k+1})} - \gamma_k^2\sqn{g_k - \nabla f(x_k)}  \\
    & & -  \gamma_k^2\sqn{\nabla f_{\xi_k}(x_{k+1}) -\nabla f(x_{k+1}) - \nabla f_{\xi_k}(x_{k}) + \nabla f(x_k)}\\
    &\geq& \frac{\gamma_k^2}{3}\sqn{\nabla f(x_{k+1})} - \gamma_k^2\sqn{g_k - \nabla f(x_k)} -\gamma_k^2\delta^2\sqn{x_{k+1} - x_k},
\end{eqnarray*}
where we used a variant of Jensen's inequality $3(a^2 + b^2 + c^2) \geq (a+b+c)^2$, for $a,b,c > 0$.
Therefore, we have 
\begin{eqnarray*}
     \sqn{x_{k+1}-x_k} &\geq& \frac{1}{1+\gamma_k^2\delta^2}\left(\frac{\gamma_k^2}{3}\sqn{\nabla f(x_{k+1})} - \gamma_k^2\sqn{g_k - \nabla f(x_k)}\right)\\
     &\geq& \frac{16}{17} \left(\frac{\gamma_k^2}{3}\sqn{\nabla f(x_{k+1})} - \gamma_k^2\sqn{g_k - \nabla f(x_k)}\right)\\
     &\geq& \frac{\gamma_k^2}{4}\sqn{\nabla f(x_{k+1})} - \gamma_k^2\sqn{g_k - \nabla f(x_k)}. 
\end{eqnarray*}
 Thus, we have 
 \begin{equation*}
     \Exp{\sqn{x_{k+1}-x_k}} \geq \frac{\gamma_k^2}{4}\Exp{\sqn{\nabla f(x_{k+1})}} - \gamma_k^2\Exp{\sqn{g_k - \nabla f(x_k)}}. 
 \end{equation*}

This concludes the proof of the lemma.

\subsection{Proof of \Cref{lemma:step3inexact}}\label{proof:step3inexact}

Let $x_{k+1} = \aprox{\epsilon}{x_k,g_k,\gamma_k,\xi_k}$.  
Then, from the definition of the function $\phi_{k}$ \eqref{eq:inexact}, we have
\begin{eqnarray*}
    \sqn{x_{k+1}-x_k} &=& \gamma_k^2 \sqn{\nabla f_{\xi_k}(x_{k+1}) + g_k - \nabla f_{\xi_k}(x_{k}){ - \nabla \phi_k(x_{k+1})}}\\
    &\geq& \gamma_k^2 \br{\frac{1}{{4}}\sqn{\nabla f(x_{k+1})} - \sqn{g_k - \nabla f(x_k)} - \delta^2\sqn{x_{k+1} - x_k} {- \sqn{\nabla \phi_k(x_{k+1})}}}.
\end{eqnarray*}
 Since $\sqn{a_1+a_2+a_3+a_4} \leq 4\left(\left\|a_1\right\|^2+\left\|a_2\right\|^2+\left\|a_3\right\|^2+\left\|a_4\right\|^2\right)$ for any vectors $a_i \in \mathbb{R}^d$, which implies $\left\|a_4\right\|^2 \geq \frac{1}{4}\left\|a_1+a_2+a_3+a_4\right\|^2-$ $\left\|a_1\right\|^2-\left\|a_2\right\|^2-\left\|a_3\right\|^2$
and  $\Exp{\sqn{\nabla \phi_k(x_{k+1})}}\leq \epsilon^2$, we deduce
\begin{equation*}
    \sqn{x_{k+1}-x_k} \geq 
\frac{\gamma_k^2}{1+\gamma_k^2\delta^2} \left(\frac{1}{4}\sqn{\nabla f(x_{k+1})} - \sqn{g_k - \nabla f(x_k)} {- \epsilon}\right).  
\end{equation*}

Therefore, we have 
\begin{eqnarray*}
     \sqn{x_{k+1}-x_k} &\geq& \frac{1}{1+\gamma_k^2\delta^2}\left(\frac{\gamma_k^2}{4}\sqn{\nabla f(x_{k+1})} - \gamma_k^2\sqn{g_k - \nabla f(x_k) }- \gamma_k^2\epsilon^2\right)\\
     &\geq& \frac{16}{17} \left(\frac{\gamma_k^2}{4}\sqn{\nabla f(x_{k+1})} - \gamma_k^2\sqn{g_k - \nabla f(x_k)}- \gamma_k^2\epsilon^2\right)\\
     &\geq& \frac{\gamma_k^2}{5}\sqn{\nabla f(x_{k+1})} - \gamma_k^2\sqn{g_k - \nabla f(x_k)} - \gamma_k^2\epsilon^2. 
\end{eqnarray*}
 Taking expectations on both sides leads to
 \begin{equation*}
     \Exp{\sqn{x_{k+1}-x_k}} \geq \frac{\gamma_k^2}{5}\Exp{\sqn{\nabla f(x_{k+1})}} - \gamma_k^2\Exp{\sqn{g_k - \nabla f(x_k)}} - \gamma_k^2\epsilon^2. 
 \end{equation*}

This concludes the proof of the lemma.

\subsection{Proof of \Cref{lemma:step_2-pp}}\label{app:proof-lemma:step_2-pp}
Recalling that 
$x^{\xi_k}_{k+1} = \aprox{\epsilon}{x_k,g_k,\gamma_k,\xi_{k}}$, we have 
\begin{equation*}
    f_{\xi_k}(x^{\xi_k}_{k+1})+\abr{g_k-\nabla f_{\xi_k} (x_k), x^{\xi_k}_{k+1} - x_k}+ \frac{1}{2\gamma_k}\|x^{\xi_k}_{k+1}-x_k\|^2 \leq f_{\xi_k}(x_{k}).
\end{equation*}

Similar to the proof of \Cref{prop:lyapunov} we start with
\begin{eqnarray*}
    f(x_{k+1}) - f(x_k) &=& \int^1_0\abr{\nabla f(\underbrace
    {x_k +\tau(x_{k+1}-x_k)}_{\eqdef x(\tau)}), x_{k+1} - x_k}d\tau,\\
    f_{\xi_k}(x^{\xi_k}_{k+1}) - f_{\xi_k}(x_k) &=& \int^1_0\abr{\nabla f_{\xi_k}(\underbrace
    {x_k +\tau(x^{\xi_k}_{k+1}-x_k)}_{\eqdef x^{\xi_k}(\tau)}), x^{\xi_k}_{k+1} - x_k}d\tau.
\end{eqnarray*}
Thus, we have 
\begin{eqnarray*}
    f(x_{k+1}) - f(x_k) &=& f_{\xi_k}(x^{\xi_k}_{k+1}) - f_{\xi_k}(x_k) \\
    &&+\int^1_0\abr{\nabla f(x(\tau)) , x_{k+1} - x_k}d\tau \\
    &&+\int^1_0\abr{-\nabla f_{\xi_k}(x^{\xi_k}(\tau)), x^{\xi_k}_{k+1} - x_k}d\tau \\
    &=& f_{\xi_k}(x^{\xi_k}_{k+1}) - f_{\xi_k}(x_k) + \abr{g_k - \nabla f_{\xi_k}(x_k), x^{\xi_k}_{k+1}-x_k}\\
    &&+\int^1_0\abr{\nabla f(x(\tau)) , x_{k+1} - x_k}d\tau \\
    &&+\int^1_0\abr{ -\nabla f_{\xi_k}(x^{\xi_k}(\tau)) - g_k + \nabla f_{\xi_k}(x_k) , x^{\xi_k}_{k+1} - x_k}d\tau.
\end{eqnarray*}
Applying the descent property of a-prox (see \Cref{def:inexact}), we deduce the following:
\begin{eqnarray*}
    f(x_{k+1}) - f(x_k) &\leq& -\frac{1}{2\gamma_k}\sqn{x^{\xi_k}_{k+1}-x_k} + \abr{\nabla f(x_k) - g_k, x^{\xi_k}_{k+1} - x_k} \\
    &&+\int^1_0\abr{\nabla f(x(\tau)) , x_{k+1} - x_k}d\tau \\
    && + \int^1_0\abr{-\nabla f_{\xi_k}(x^{\xi_k}(\tau)) - \nabla f(x_k) + \nabla f_{\xi_k}(x_k) , x^{\xi_k}_{k+1} - x_k}d\tau.
\end{eqnarray*}

Let us take  expectation from both sides  conditioned to $\cG_k = \curlybr{x_k,x_{k+1},S_k,g_k}$. In other words, we take expectation with respect to the random index $\xi_k$ chosen uniformly  from the already chosen $S_k$:
\begin{eqnarray*}
    f(x_{k+1}) - f(x_k)
    &\leq& \Exp{-\frac{1}{2\gamma_k}\sqn{x^{\xi_k}_{k+1}-x_k} + \abr{\nabla f(x_k) - g_k, x^{\xi_k}_{k+1} - x_k}\mid \cG_k} \\
    &&+ \Exp{\int^1_0\abr{\nabla f(x(\tau)) , x_{k+1} - x_k}d\tau \mid \cG_k} \\
    && + \ExpCond{\int^1_0\abr{-\nabla f_{\xi_k}(x^{\xi_k}(\tau)) - \nabla f(x_k) + \nabla f_{\xi_k}(x_k) , x^{\xi_k}_{k+1} - x_k}d\tau }{ \cG_k} \\
    &=& \Exp{-\frac{1}{2\gamma_k}\sqn{x^{\xi_k}_{k+1}-x_k}\mid \cG_k} + \abr{\nabla f(x_k) - g_k, x_{k+1} - x_k} \\
    && + \Exp{\int^1_0\abr{\nabla f(x(\tau)) - \nabla f(x_k) -\nabla f_{\xi_k}(x^{\xi_k}(\tau)) + \nabla f_{\xi_k}(x_k) , x^{\xi_k}_{k+1} - x_k}d\tau \mid \cG_k}.
\end{eqnarray*}
Here the last equality is due to the fact that $\xi_k$ is independent of $\cG_k$ and $x_{k+1} = \Exp{x^{\xi_k}_{k+1} \mid \cG_k}$.
Therefore, applying the Cauchy-Schwartz inequality
\begin{align*}
    f(x_{k+1}) &-f(x_k) \leq   \Exp{-\frac{1}{2\gamma_k}\sqn{x^{\xi_k}_{k+1}-x_k}\mid \cG_k} + \abr{\nabla f(x_k) - g_k, x_{k+1} - x_k} & \\
    & \quad + \Exp{\int^1_0\abr{\nabla f(x(\tau)) - \nabla f(x^{\xi_k}(\tau)), x^{\xi_k}_{k+1} - x_k}d\tau \mid \cG_k} \\
    & \quad + \Exp{\int^1_0\abr{\nabla f(x^{\xi_k}(\tau)) - \nabla f(x_k) -\nabla f_{\xi_k}(x^{\xi_k}(\tau)) + \nabla f_{\xi_k}(x_k) , x^{\xi_k}_{k+1} - x_k}d\tau \mid \cG_k}\\
    &\leq  \Exp{-\frac{1}{2\gamma_k}\sqn{x^{\xi_k}_{k+1}-x_k}\mid \cG_k} + \abr{\nabla f(x_k) - g_k, x_{k+1} - x_k} \\
    & \quad + \Exp{\int^1_0\norm{\nabla f(x(\tau)) - \nabla f(x^{\xi_k}(\tau))}\norm{x^{\xi_k}_{k+1} - x_k}d\tau \mid \cG_k} \\
    & \quad + \Exp{\int^1_0\norm{\nabla f(x^{\xi_k}(\tau)) - \nabla f(x_k) -\nabla f_{\xi_k}(x^{\xi_k}(\tau)) + \nabla f_{\xi_k}(x_k) }\norm{ x^{\xi_k}_{k+1} - x_k}d\tau \mid \cG_k}.
\end{align*}
Applying Cauchy-Schwartz inequality once again, we deduce
\begin{eqnarray*}
    f(x_{k+1}) - f(x_k)
    &\leq & \Exp{-\frac{1}{2\gamma_k}\sqn{x^{\xi_k}_{k+1}-x_k}\mid \cG_k} + \frac{C}{2}\sqn{\nabla f(x_k) - g_k} + \frac{1}{2C}\sqn{x_{k+1} - x_k}\\
    && + \Exp{\int^1_0 L \norm{x(\tau) - x^{\xi_k}(\tau)}\norm{x^{\xi_k}_{k+1} - x_k}d\tau \mid \cG_k} \\
    && + \Exp{\int^1_0 \delta \norm{x^{\xi_k}(\tau) - x_k}\norm{ x^{\xi_k}_{k+1} - x_k}d\tau \mid \cG_k}\\
    &\leq & \Exp{-\frac{1}{2\gamma_k}\sqn{x^{\xi_k}_{k+1}-x_k}\mid \cG_k} + \frac{C}{2}\sqn{\nabla f(x_k) - g_k} + \frac{1}{2C}\sqn{x_{k+1} - x_k}\\
    && + \Exp{\int^1_0 L \tau\norm{x_{k+1} - x^{\xi_k}_{k+1}}\norm{x^{\xi_k}_{k+1} - x_k}d\tau \mid \cG_k} \\
    && + \Exp{\int^1_0 \delta \tau\normsq{ x^{\xi_k}_{k+1} - x_k}d\tau \mid \cG_k}.
\end{eqnarray*}
Computing the integral with respect to $\tau$, we obtain
\begin{eqnarray*}
    f(x_{k+1}) - f(x_k)
    &\leq & \Exp{-\frac{1}{2\gamma_k}\sqn{x^{\xi_k}_{k+1}-x_k}\mid \cG_k} + \frac{C}{2}\sqn{\nabla f(x_k) - g_k} + \frac{1}{2C}\sqn{x_{k+1} - x_k}\\
    && + \frac{L}{2} \Exp{\norm{x_{k+1} - x^{\xi_k}_{k+1}}\norm{x^{\xi_k}_{k+1} - x_k}\mid \cG_k}  + \frac{\delta }{2} \Exp{ \normsq{ x^{\xi_k}_{k+1} - x_k} \mid \cG_k}\\
     &\leq & \Exp{-\frac{1}{2\gamma_k}\sqn{x^{\xi_k}_{k+1}-x_k}\mid \cG_k} + \frac{C}{2}\sqn{\nabla f(x_k) - g_k} + \frac{1}{2C}\sqn{x_{k+1} - x_k}\\
    && + \frac{L}{4} \Exp{\normsq{x_{k+1} - x^{\xi_k}_{k+1}}\mid \cG_k}  + \frac{2\delta + L}{4} \Exp{ \normsq{ x^{\xi_k}_{k+1} - x_k} \mid \cG_k}.
\end{eqnarray*}
Recall again that $x_{k+1} = \Exp{x^{\xi_k}_{k+1} \mid \cG_k}$, thus $x_{k+1} = \argmin_{a\in \R^d} \Exp{\normsq{x^{\xi_k}_{k+1} - a } \mid \cG_k}$. Therefore, 
\begin{equation*}
    \Exp{\normsq{x^{\xi_k}_{k+1} - x_{k+1}}\mid \cG_k} 
    \leq \Exp{\normsq{x^{\xi_k}_{k+1} - x_{k}}\mid \cG_k}.
\end{equation*}
Furthermore,
\begin{equation*}
    \sqn{x_{k+1} - x_k} = \sqn{\Exp{x^{\xi_k}_{k+1} \mid \cG_k} - x_k} 
    \leq \Exp{\sqn{x^{\xi_k}_{k+1} - x_k}  \mid \cG_k}.
\end{equation*}
Combining these two bounds, we deduce
\begin{eqnarray*}
    f(x_{k+1}) - f(x_k)
     &\leq &  \frac{C}{2}\sqn{\nabla f(x_k) - g_k} \\
    &&   + \roundbr{\frac{1}{2C} + \frac{\delta + L}{2} -\frac{1}{2\gamma_k}} \Exp{ \normsq{ x^{\xi_k}_{k+1} - x_k} \mid \cG_k}.
\end{eqnarray*}
The previous bound is true for every positive value of $C$. Thus, it is true also for $C = 4\gamma_k$. Taking into account that $\gamma_k < \frac{1}{4(L + \delta)}$, we get 
\begin{equation*}
    \frac{1}{2C} + \frac{\delta + L}{2} -\frac{1}{2\gamma_k} 
    \leq \frac{1}{8\gamma_k} + \frac{1}{8\gamma_k} -\frac{1}{2\gamma_k} 
    =  -\frac{1}{4\gamma_k}.
\end{equation*}
Therefore, 
\begin{eqnarray*}
    f(x_{k+1}) - f(x_k)
     &\leq &  {2\gamma_k}\sqn{\nabla f(x_k) - g_k}
       - \frac{1}{4\gamma_k} \Exp{ \normsq{ x^{\xi_k}_{k+1} - x_k} \mid \cG_k}.
\end{eqnarray*}

Thus, taking full expectation on both sides, we have 
\begin{equation*}
    \Exp{f(x_{k+1}) - f_{\inf}} \leq \Exp{f(x_{k}) - f_{\inf}} 
    + {2\gamma_k}\Exp{\sqn{\nabla f(x_k) - g_k}}
       - \frac{1}{4\gamma_k} \Exp{ \normsq{ x^{\xi_k}_{k+1} - x_k}}.
\end{equation*}
This concludes the proof.

\subsection{Proof of \Cref{lemma:step_3-pp}}\label{app:proof-lemma:step_3-pp}

By the definition of $x^{\xi_k}_{k+1}$, for every $\xi \in S_k$ we have 
\begin{eqnarray*}
    \sqn{x^{\xi_k}_{k+1}-x_k} 
    &=& \gamma_k^2 \sqn{\nabla f_{\xi_k}(x^{\xi_k}_{k+1}) + g_k - \nabla f_{\xi_k}(x_{k})  - \nabla \phi_{k} (x_{k+1})}\\
    &=& \gamma_k^2 \left\| \nabla f(x^{\xi_k}_{k+1}) +g_k - \nabla f(x_k) + \nabla f_{\xi_k}(x^{\xi_k}_{k+1}) 
    \right. \\ 
    & & \left. - \nabla f(x^{\xi_k}_{k+1}) - \nabla f_{\xi_k}(x_{k}) + \nabla f(x_k) - \nabla \phi_{k} (x_{k+1})\right\|^2\\
    &\geq&\frac{\gamma_k^2}{4}\sqn{\nabla f(x^{\xi_k}_{k+1})} - \gamma_k^2\sqn{g_k - \nabla f(x_k)}  \\ 
    && -\gamma_k^2\sqn{\nabla f_{\xi_k}(x^{\xi_k}_{k+1}) -\nabla f(x^{\xi_k}_{k+1}) - \nabla f_{\xi_k}(x_{k}) + \nabla f(x_k)} - \gamma_k^2\epsilon^2\\
    &\geq& \frac{\gamma_k^2}{4}\sqn{\nabla f(x^{\xi_k}_{k+1})} - \gamma_k^2\sqn{g_k - \nabla f(x_k)} -\gamma_k^2\delta^2\sqn{x^{\xi_k}_{k+1} - x_k}- \gamma_k^2\epsilon^2.
\end{eqnarray*}
The third inequality is due to Cauchy-Schwartz and the second property of the approximate proximal operator (See \Cref{def:inexact}).
Therefore, we have 
\begin{eqnarray*}
     \sqn{x^{\xi_k}_{k+1}-x_k} &\geq& \frac{1}{1+\gamma_k^2\delta^2}\left(\frac{\gamma_k^2}{4}\sqn{\nabla f(x^{\xi_k}_{k+1})} - \gamma_k^2\sqn{g_k - \nabla f(x_k)}- \gamma_k^2\epsilon^2\right)\\
     &\geq& \frac{16}{17} \left(\frac{\gamma_k^2}{4}\sqn{\nabla f(x^{\xi_k}_{k+1})} - \gamma_k^2\sqn{g_k - \nabla f(x_k)}- \gamma_k^2\epsilon^2\right)\\
     &\geq& \frac{\gamma_k^2}{5}\sqn{\nabla f(x^{\xi_k}_{k+1})} - \gamma_k^2\sqn{g_k - \nabla f(x_k)}- \gamma_k^2\epsilon^2. 
\end{eqnarray*}
We deduce 
 \begin{equation*}
     \Exp{\sqn{x^{\xi_k}_{k+1}-x_k}} \geq \frac{\gamma_k^2}{5}\Exp{\sqn{\nabla f(x^{\xi_k}_{k+1})}} - \gamma_k^2\Exp{\sqn{g_k - \nabla f(x_k)}} - \gamma_k^2\epsilon^2. 
 \end{equation*}
This concludes the proof of the lemma. 

\section{Experimental details} \label{app:experiment}
We provide additional details on the experimental settings from Section \ref{sec:experiments}.

Consider a distributed ridge regression problem defined as
\begin{equation} \label{eq:quad_problem}
    f(x) = \ExpSub{\xi}{\sqn{A_{\xi}x - y_{\xi}}} + \frac{\lambda}{2} \sqn{x},
\end{equation}
where $\xi$ is uniform random variable over $\{1, \dots, n\}$ for $n = 10, \lambda=0.1$.
We follow a similar to \cite{lin2024stochastic} procedure for synthetic data generation, which allows us to calculate and control Hessian similarity $\delta$.
Namely, a random matrix $A_0 \in \mathbb{R}^{d\times d}$ ($d=100$) is generated with entries from a standard Gaussian distribution $\mathcal{N}(0, 1)$. Then we obtain $A = A_0 A_0^{\top}$ (to ensure symmetry) and set $A'_{\xi} = A + B_{\xi}$ by adding a random symmetric matrix $B_{\xi}$ (generated similarly to $A$). Afterwards we modify $A_{\xi} = A'_{\xi} + I \lambda_{\min} (A'_{\xi})$ by adding an identity matrix $I$ times minimum eigenvalue to guarantee $A_{\xi} \succeq 0$. Entries of vectors $y_\xi \in \sR^d$, and initialization $x_0 \in \sR^d$ are generated from a standard Gaussian distribution $\mathcal{N}(0, 1)$.

In the case of inexact proximal point computation (1/10 local steps), local subproblems \eqref{eq:inexact} are solved with gradient descent.

Simulations were performed on a machine with $24 \operatorname{Intel}(\mathrm{R}) \operatorname{Xeon}(\mathrm{R})$ Gold 6246 $\mathrm{CPU}$ @ 3.30 $\mathrm{GHz}$.

\end{document}